\documentclass[11pt, reqno]{amsart}
\usepackage{mathrsfs}
\usepackage{lmodern}
\usepackage[lmargin=1in,rmargin=1in,tmargin=1in,bmargin=1in,marginparwidth=5 in, marginparsep= .2 in]{geometry}
\usepackage[dvipsnames]{xcolor}
\usepackage[utf8]{inputenc}
\usepackage{bbold}
\usepackage[colorlinks=true,
            linkcolor=blue,
            citecolor=blue,
            urlcolor=cyan]{hyperref}
\usepackage{comment}
\usepackage{amsmath,amsfonts,blindtext, amssymb,amsbsy,amscd,amsthm,bbm,enumerate,epsfig,esint,graphicx,float,marginnote,mathrsfs,mathtools,multicol}
\usepackage[shortlabels]{enumitem}

\usepackage{tikz}
\usepackage{microtype}

\usepackage{thm-restate}

\usepackage{parskip}
\allowdisplaybreaks

 \newcommand{\bC}{\mathbb{C}}

\newcommand{\bR}{\mathbb{R}}\newcommand{\bZ}{\mathbb{Z}}\newcommand{\cA}{\mathcal{A}}\newcommand{\cB}{\mathcal{B}}\newcommand{\cD}{\mathcal{D}}\newcommand{\cE}{\mathcal{E}}\newcommand{\cF}{\mathcal{F}}\newcommand{\cN}{\mathcal{N}}
 \newcommand{\SL}{\mathrm{SL}}
\providecommand{\cW}{\mathcal{W}}
\DeclareMathOperator{\Cov}{Cov}
\newcommand{\gen}{\mathrm{gen}}
\renewcommand{\Re}{\operatorname{Re}}

\makeatletter
\renewenvironment{proof}[1][\proofname]{\par
  \vspace{-.7\topsep}% remove the space after the theorem
  \pushQED{\qed}%
  \normalfont
  \topsep0pt \partopsep0pt % no space before
  \trivlist
  \item[\hskip\labelsep
        \itshape
    #1\@addpunct{.}]\ignorespaces
}{%
  \popQED\endtrivlist\@endpefalse
  \addvspace{6pt plus 6pt} % some space after
}
\makeatother

\newtheorem{theorem}{Theorem}[section]
\newtheorem{prop}[theorem]{Proposition}
\newtheorem{cor}[theorem]{Corollary}
\newtheorem{lemma}[theorem]{Lemma}
\newtheorem{definition}[theorem]{Definition}
\newtheorem{example}[theorem]{Example}

\newtheorem{conjecture}[theorem]{Conjecture}
\newtheorem{remark}[theorem]{Remark}

\newcommand{\nc}{\newcommand}
\nc{\dmo}{\DeclareMathOperator}

\let\oldtocsection=\tocsection

\let\oldtocsubsection=\tocsubsection

\let\oldtocsubsubsection=\tocsubsubsection

\renewcommand{\tocsection}[2]{\hspace{0em}\oldtocsection{#1}{#2}}
\renewcommand{\tocsubsection}[2]{\hspace{2em}\oldtocsubsection{#1}{#2}}
\renewcommand{\tocsubsubsection}[2]{\hspace{4.5em}\oldtocsubsubsection{#1}{#2}}
\usepackage{xpatch}
\makeatletter   
\xpatchcmd{\@tocline}
{\hfil\hbox to\@pnumwidth{\@tocpagenum{#7}}\par}
{\ifnum#1<0\hfill\else\dotfill\fi\hbox to\@pnumwidth{\@tocpagenum{#7}}\par}
{}{}
\makeatother

\makeatletter
\renewcommand*\@maketitle{%
  \normalfont\normalsize
  \@adminfootnotes
  \@mkboth{\@nx\shortauthors}{\@nx\shorttitle}%
  \global\topskip42\p@\relax % 5.5pc   "   "   "     "     "
  \@settitle
  \ifx\@empty\authors \else \@setauthors \fi
  \ifx\@empty\@dedicatory
  \else
    \baselineskip18\p@
    \vtop{\centering{\footnotesize\itshape\@dedicatory\@@par}%
      \global\dimen@i\prevdepth}\prevdepth\dimen@i
  \fi
  \@setabstract
  \normalsize
  \if@titlepage
    \newpage
  \else
    \dimen@34\p@ \advance\dimen@-\baselineskip
    \vskip\dimen@\relax
  \fi
}
\makeatother

\nc{\C}{\mathbb{C}}
\nc{\F}{\mathbb{F}}
\nc{\N}{\mathbb{N}}
\nc{\Q}{\mathbb{Q}}
\nc{\R}{\mathbb{R}}
\nc{\T}{\mathbb{T}}
\nc{\Z}{\mathbb{Z}}
\newcommand{\eps}{\varepsilon}

\DeclareMathOperator{\dist}{dist}

\DeclareMathAlphabet{\mathcal}{OMS}{cmsy}{m}{n}
\numberwithin{equation}{section}

\renewcommand{\hat}{\widehat}
\renewcommand{\tilde}{\widetilde}

\renewcommand{\bar}{\overline}

\newcommand{\del}{\partial}

\newcommand{\abs}[1]{\left\lvert #1 \right\rvert}
\newcommand{\norm}[1]{\lVert #1 \rVert}

\newcommand\angles[1]{\left \langle #1 \right \rangle}
\newcommand\brax[1]{\left \{ #1 \right \} }
\newcommand\floor[1]{\left \lfloor #1 \right \rfloor }
\newcommand\paren[1]{\left ( #1 \right ) }
\newcommand\braks[1]{\left [ #1 \right ] }

\title[Breaking Universality in the Lower Order Terms in the 1-level and 2-level Density of Holomorphic Cusp Newforms]{Breaking Universality in the Lower Order Terms in the 1-level and 2-level Density of Holomorphic Cusp Newforms}

\date{\today}

\author{Yuval Amit}
\address{Department of Mathematics, University of Maryland, College Park, MD 20742}
\email{\href{mailto:}{yamit@umd.edu}}

\author{Lucas Chen}
\address{Department of Mathematics, University of Chicago, Chicago, IL 60637}
\email{\href{mailto:}{lucasch@uchicago.edu}}

\author{Lawrence Dillon}
\address{Department of Mathematics, University of Washington, Seattle, WA 98105}
\email{\href{mailto:}{lvid@uw.edu}}

\author{Christopher Housholder}
\address{Department of Mathematics, Missouri State University, Springfield, MO 65897}
\email{\href{mailto:}{christopherlancehousholder@gmail.com}}

\author{Xiaoyao Huang}
\address{Department of Mathematics, University of Michigan, Ann Arbor, MI 48109}
\email{\href{mailto:}{xyrushac@umich.edu}}

\author{Joshua Khan}
\address{Department of Mathematics, University of Rochester, Rochester, NY 14627}
\email{\href{mailto:}{joshuakhan1307@gmail.com}}

\author{Say-Yeon Kwon}
\address{Department of Mathematics, Princeton University, Princeton, NJ 08544}
\email{\href{mailto:sk9017@princeton.edu}{sk9017@princeton.edu}}

\author{Meiling Laurence}
\address{Department of Mathematics, Yale University, New Haven, CT 06520}
\email{\href{mailto:}{meiling.laurence@yale.edu}}

\author{Steven J. Miller}
\address{Department of Mathematics \& Statistics, Williams College, Williamstown, MA 01267}
\email{\href{mailto:sjm1@williams.edu}{sjm1@williams.edu}}

\author{Vishal Muthuvel}
\address{Department of Mathematics, Columbia University, New York, NY 10027}
\email{\href{mailto:vm2696@columbia.edu}{vm2696@columbia.edu}}

\author{Devayani Pradhan}
\address{Department of Mathematics, University of Michigan, Ann Arbor, MI 48109}
\email{\href{mailto:}{pradhand@umich.edu}}

\author{Luke Rowen}
\address{Department of Mathematics, Carleton College, Northfield, MN 55057 }
\email{\href{mailto:rowenl@carleton.edu}{rowenl@carleton.edu}}

\author{Pramana Saldin}
\address{Department of Mathematics, University of Wisconsin, Madison, WI 53706}
\email{\href{mailto:saldin@wisc.edu}{saldin@wisc.edu}}

\author{Connor Yau}
\address{Department of Mathematics \& Statistics, Williams College, Williamstown, MA 01267}
\email{\href{mailto:cy10@williams.edu}{cy10@williams.edu}}

\author{Steven Zanetti}
\address{Department of Mathematics, University of Michigan, Ann Arbor, MI 48109}
\email{\href{mailto:szanetti@umich.edu}{szanetti@umich.edu}}

\begin{document}

%%%%%%%%%%%%%%%%

\begin{abstract}
The Katz--Sarnak density conjecture predicts that the normalized low-lying zeros of a family of
$L$-functions ordered by analytic conductor $R$ converge, as $R \to \infty$, to the eigenvalue
statistics near $1$ of a classical compact group, and that the limit depends on the family only
through its symmetry type. We study the first correction to this limit for families of holomorphic
cusp newforms of fixed weight $k$ and varying level $N$. Sharpening work of Miller on lower order
terms in the one-level density, we evaluate the harmonically weighted one- and two-level densities
with errors $O\paren{\log^{-5}(R)}$ and $O\paren{\log^{-2}(R)}$, and show that each of them already
differs at order $\log^{-1}(R)$ from an explicit universal profile built out of $k$, the test
functions, and $R$ alone. We compute that difference exactly. It is the product of a functional of
the test functions with a single arithmetic invariant $\kappa_1(\cF)$ of the family, an explicit sum
over the bounded prime factors of the level. Two families of newforms of the same weight and the
same symmetry type, with conductors of the same size, therefore have densities whose first lower
order terms differ by an explicit nonzero amount whenever their levels have different bounded prime
factors.
\end{abstract}

\maketitle
% \vspace{-13mm}
\tableofcontents
\pagestyle{plain}

%%%%%%%%%%%%%%%%%%%%%%%%%%%%%%%%%%%%%%%%%%%%%%%%%%%%%%%%%%
%%%%%%%%%%%%%%%%%%%%%%%%%%%%%%%%%%%%%%%%%%%%%%%%%%%%%%%%%%
%%%%%%%%%%%%%%%%%%%%%%%%%%%%%%%%%%%%%%%%%%%%%%%%%%%%%%%%%%
\section{Introduction}

Since the 1970s, the study of zeros of $L$-functions has revealed striking connections with random matrix theory. In particular, the local statistics of zeros high on the critical line are conjecturally governed by the eigenvalue statistics of the Gaussian Unitary Ensemble, a phenomenon often referred to as universality. The behavior of zeros near the central point $s=1/2$, however, is more sensitive to the arithmetic of the underlying family. This distinction is particularly important for questions such as the Birch and Swinnerton--Dyer conjecture, which are governed by the order of vanishing at the central point \cite{Hej, Mon, Od1, Od2, RS}.

The Katz--Sarnak philosophy predicts that, for a family of $L$-functions ordered by analytic conductor, the statistics of zeros near the central point are governed by the eigenvalue statistics of a classical compact group determined by the family \cite{katz1999random,katz1999zeros}. This prediction has been supported by a large body of work, including results on moments \cite{CF2000,CFKRS2005,KeSn1,KeSn2,KeSn3} and on $n$-level densities for suitable test functions \cite{DM1,FI,Gu,HR,HM,Mil2,OS,RR1,Ro,Rub,Yo2}.

% Limiting behavior in general is determined by leading terms, while the lower-order terms impact the rate of convergence, and retain finer arithmetic information which is invisible in the limiting random-matrix model. Katz--Sarnak conjectures that the density $D_n(\mathcal F)$ has a universal behavior across all families of $L$-functions $\cF$. However, we show that the limiting behavior of lower-order terms is not universal for different families of $L$-functions derived from holomorphic cusp newforms. Non-universal behavior of the lower order terms can reveal

Given an $L$-function $L(s,f)$ associated with a holomorphic cusp newform $f \in H^*_k(N) $ of level $N$ and weight $k$, we assume the Generalized Riemann Hypothesis (GRH).\footnote{We do not need the GRH for any of the calculations. However, assuming the GRH allows us to order the zeros and relate the $n$-level density with the spacing between the zeros.} Thus, all non-trivial zeros of $L(s,f)$ are of the form $\rho_f=1/2+i\gamma_f$. We then enumerate zeros by their imaginary part $0 \leq \gamma_f^{(1)}\leq \gamma_f^{(2)}\leq \cdots$ with $\gamma_f^{(-j)}=-\gamma_f^{(j)}$ by symmetry. We define the $n$-level density by 
 \begin{equation}\label{eq:intro-1}
    D_n(f;\Phi)\ :=\ \sum_{\substack{j_1,\dots,j_n\\j_i \neq \pm j_k}}\phi_1\left( \gamma^{(j_1)}_f \frac{\log (R)}{2\pi}\right)\cdots \phi_n\left( \gamma^{(j_n)}_f \frac{\log (R)}{2\pi}\right),
 \end{equation}
 where $\Phi(x_1, \dots, x_n)=\phi_1(x_1)\cdot \phi_2(x_2)\cdots\phi_n(x_n)$ is a Schwartz test function\footnote{In general, we do not need the test function $\Phi$ to be a product of coordinate-wise test functions \cite{cheek2026densitylowlyingzeros}.} and $R$ is the analytic conductor of $f$. We rescale the zeros near the central point by $\log (R)/2\pi$. Note that for $f \in H^{*}_k(N)$, its analytic conductor $R$ is $(64\pi^2)^{-1}N(k+1)(k+3)$, and so $\log(R) \sim\log(N)$. 
\begin{conjecture}[Katz-Sarnak \cite{katz1999random}]
    For a family $\mathcal{F}=\bigcup\mathcal{F}_N$ of $L$-functions ordered by their conductors, we have
 \begin{equation}\label{eq:intro-2}
     \lim_{N\to \infty}\frac{\sum_{f\in \mathcal{F}_N}D_n(f;\Phi)}{|\mathcal{F}_N|}\ =\ \int \cdots\int \Phi(x_1,\dots, x_n)W_{n, G(\mathcal{F})}(x_1,\dots,x_n)dx_1\cdots dx_n,
 \end{equation}
where $W_{G(\mathcal{F})}$ represents the limiting distribution of a similar statistic for the eigenvalues of random matrices in some classical compact group as their size goes to infinity.
\end{conjecture}

The Katz--Sarnak prediction concerns only the leading term of the density. Lower order terms retain finer arithmetic information that is invisible in the limiting random matrix model, and they control the rate at which the conjectural limit is approached. Families with the same symmetry type can therefore behave quite differently below the leading order, and this has been observed for several families of $L$-functions. For instance, Young \cite{young2005lower} showed that the lower order terms in families of elliptic curves depend strongly on the arithmetic of the family. See also \cite{GAO_ZHAO_2023, Fiorilli_2015, FI} for lower order terms in other families.

We study the lower order behavior of the 1- and 2-level densities of zeros for families of $L$-functions attached to holomorphic cusp newforms. Taking the $L$-functions of level-$N$ weight-$k$ newforms, we show that the lower order terms depend explicitly on the prime factorization of $N$ as $N \to \infty$.

A helpful analogy, due to Miller \cite{Miller2009_LOTerms_1LevelDensity}, comes from the Central Limit Theorem. Let $X_1,\dots,X_n$ be independent, identically distributed random variables with mean $\mu$, variance $\sigma^2$ and finite third moment. Then $(X_1+\cdots+X_n-n\mu)/(\sigma\sqrt n)$ converges to the standard normal $\cN(0,1)$ as $n \to \infty$, and the universality is that the properly normalized main term does not remember the initial distribution. The rate of convergence, on the other hand, depends on the higher moments of the $X_i$ through the Berry--Esseen theorem \cite{Berry1941, Esseen1942}. Similarly, our results study the lower order terms which are forgotten by the universal terms.

Such differences are worth chasing because they carry information about the family that the leading term cannot see. Miller \cite{Miller2009_LOTerms_1LevelDensity} found elliptic curve families whose density statistics stray from the conjectural predictions, and in \cite{Mil2, Miller2005, Miller2009_LOTerms_1LevelDensity} observed that different one-parameter families of elliptic curve $L$-functions have different lower order terms in their moments. That difference has consequences for the distribution of zeros near the central point, and it accounts for features such as excess rank.

%%%%%%%%%%%%%%%%%%%%%%%%%%%%%%%%%%%%%%%%%%%%%%%%%%%%%%%%%%%%%%%%%%%%%%
\subsection{Main results}\label{subsec:mainresults}

\begin{definition}\label{def:family}
    Let $H_k^*(N)$ be the set of cuspidal newforms of level $N$ and weight $k$. We fix $k$, and let $\cF_{N}$ be the associated set of $L$-functions defined through \eqref{eq:lfndefn}. To each $f \in \cF_{N}$ we assign a weight $w_R(f)$ defined in \eqref{eq:wR_defn}, where $R$ is the analytic conductor of $f$ (note that $R$ is the same for all $f \in \cF_{N}$). Let
\begin{align}\label{eq:WR_defn}
    W_R(\cF_{N})\ :=\  \sum_{f\in \cF_{N}}w_R(f),
\end{align}
and let $\phi,\phi_1,\phi_2$ be test functions which are even and Schwartz, and whose Fourier transforms
\begin{align}
    \hat \phi(y)\ :=\ \int_{-\infty}^\infty \phi(x) e^{-2\pi i xy}\,dx
\label{eq:family-1}
\end{align}
have compact support. We write $\langle\,\cdot\,\rangle_N$ for the harmonic average over $\cF_N$, so that for any $X \colon \cF_N \to \bC$,
\begin{align}\label{eq:avg-notation}
    \angles{X}_N\ :=\ \frac{1}{W_R(\cF_{N})}\sum_{f \in \cF_{N}}w_R(f)X(f).
\end{align}
The weighted first- and second-level densities at level $N$ are
\begin{equation}\label{eq:1stlvl-density-def}
     D_1(\cF_{N},\phi)\ :=\ \angles{D_1(\,\cdot\,;\phi)}_N,
     \qquad
     D_2(\cF_{N},\phi_1,\phi_2)\ :=\ \angles{D_2(\,\cdot\,;\phi_1,\phi_2)}_N,
\end{equation}
where
\begin{equation}\label{eq:family-2}
    D_1(f; \phi) \ :=\ \displaystyle\sum_{i}\phi\left(\frac{\log (R)}{2\pi}\gamma_f^{(i)}\right),
    \qquad
D_2(f; \phi_1,\phi_2)\ :=\ \sum_{\substack{i,j\\i \neq \pm j}}\phi_1 \left( \frac{\log (R)}{2\pi} \gamma_f^{(i)} \right)\phi_2 \left( \frac{\log (R)}{2\pi} \gamma_f^{(j)} \right).
\end{equation}
\end{definition}

Both $D_1(\cF_N,\phi)$ and $D_2(\cF_N,\phi_1,\phi_2)$ converge, as $N \to \infty$, to the Katz--Sarnak predictions for the orthogonal symmetry group, recorded in \eqref{eq:D1-profile-value} and \eqref{eq:D2-profile-value}. We study is the rate at which they converge, and the extent to which that rate is determined by the arithmetic of the level.

\begin{definition}[Factorization profile]\label{def:profile}
    Fix an integer $n \geq 1$, exponents $a_1,\dots,a_n \in \bZ^+$, and reals $\delta_1,\dots,\delta_n \geq 0$ with $\sum_i a_i\delta_i = 1$. We call an integer sequence $\{N\} \to \infty$ \textbf{admissible} if $N$ is of the form
    \begin{align}\label{eq:level-shape}
        N\ =\ q_1^{a_1}\cdots q_n^{a_n}, \qquad q_i \text{ distinct primes with } q_i \asymp N^{\delta_i}.
    \end{align}
    The $q_i$ with $\delta_i = 0$ are bounded, so after passing to a subsequence we may assume that $\delta_i=0$ implies $q_i$ is a fixed prime, independent of $N$. The associated \textbf{family} is $\cF := \bigcup_N \cF_{N}$, and we set
    \begin{align}\label{eq:Q-defs}
        Q_0(\cF)\ :=\ \brax{q_i: \delta_i = 0}, \qquad
        Q(\cF)\ :=\ \brax{q_i : \delta_i = 0 \text{ and } a_i = 1}\ \subseteq\ Q_0(\cF).
    \end{align}
    Thus $Q_0(\cF)$ is the set of bounded primes dividing the level, and $Q(\cF)$ the subset of those dividing it exactly once. Both are finite and fixed along the sequence.
\end{definition}

When analyzing the densities for admissible families of newforms, we break down our expressions of interest into computable invariants whose values explicitly depend on the prime factorization of the level $N$.

\begin{definition}[Level invariants]\label{def:kappas}
    For an admissible family $\cF$ with $Q_0 = Q_0(\cF)$ and $Q = Q(\cF)$, set
    \begin{align}
        \kappa_1(\cF)\ &:=\ \sum_{q \in Q_0}\frac{\log (q)}{q}\ +\ \sum_{q \in Q}\frac{\log (q)}{q^2-1}, \label{eq:kappa1-def}\\
        \kappa_2(\cF)\ &:=\ \sum_{q \in Q_0}\frac{(q-1)\log^2 (q)}{q^2}\ -\ \sum_{q \in Q}\frac{q^2\log^2( q)}{\paren{q^2-1}^2}, \label{eq:kappa2-def}\\
        \kappa_3(\cF)\ &:=\ \sum_{q \in Q_0}\frac{\log^3 (q)}{q}\ +\ \sum_{q \in Q}\frac{q^2\paren{q^2+1}\log^3 (q)}{\paren{q^2-1}^3}. \label{eq:kappa3-def}
    \end{align}
\end{definition}

Every summand of $\kappa_1(\cF)$ is positive, so
\begin{align}\label{eq:kappa1-positive}
    \kappa_1(\cF)\ >\ 0 \iff Q_0(\cF) \neq \emptyset \iff \text{the level has a bounded prime factor.}
\end{align}

Our results compare the densities of a family against an explicit universal profile, by which we mean a function of the weight $k$, the test functions, and the analytic conductor $R$ alone. We define $\cD_1(R;\phi)$ and $\cD_2(R;\phi_1,\phi_2)$ in \eqref{eq:D1-profile-def} and \eqref{eq:D2-profile-def}, which resemble the predicted behavior according to the Katz--Sarnak conjecture. The theorems below show that the deviation from the profile is nonzero, that it is visible already at the first lower order term, and that it is an explicit function of the bounded part of the factorization of $N$.

\begin{restatable}{theorem}{FirstUniversalityBreaks}\label{thm:D1-universality-breaks}
    Assume the Generalized Riemann Hypothesis for the $L$-functions in $\cF$, which we use only to order the zeros. Let $\cF$ be an admissible family as in Definition \ref{def:profile}, and let $\phi$ be an even Schwartz function with $\hat\phi$ supported in $[-\sigma,\sigma]$ for some $\sigma < 4/3$. Then
    \begin{align}\label{eq:D1-master}
        D_1(\cF_{N},\phi)\ =\ \cD_1(R;\phi)\ -\ \frac{2\hat\phi(0)\kappa_1(\cF)}{\log (R)}\ -\ \frac{4\hat\phi''(0)\kappa_3(\cF)}{\log^3 (R)}\ +\ O\paren{\log^{-5}(R)}.
    \end{align}
    In particular
    \begin{align}\label{eq:D1-invariant}
        \lim_{N \to \infty}\log (R)\braks{D_1(\cF_{N},\phi)-\cD_1(R;\phi)}\ =\ -2\hat\phi(0)\kappa_1(\cF),
    \end{align}
    which is nonzero whenever $\hat\phi(0) \neq 0$ and the level has a bounded prime factor. Consequently, if $\cF$ and $\tilde\cF$ are two admissible families whose levels run through sequences with $N \asymp \tilde N$, then
    \begin{align}\label{eq:D1-comparison}
        \lim_{N \to \infty} \log (R)\braks{D_1(\cF_{N},\phi)-D_1(\tilde\cF_{\tilde N},\phi)}\ =\ -2\hat\phi(0)\paren{\kappa_1(\cF)-\kappa_1(\tilde\cF)}.
    \end{align}
\end{restatable}

Since $\hat\phi(0) = \int\phi$, the hypothesis $\hat\phi(0) \neq 0$ holds for every nonnegative $\phi \not\equiv 0$, and in particular for the standard test functions used in this subject.

The second level is governed by the same invariant $\kappa_1(\cF)$, paired with a different functional of the test functions. For even Schwartz $\phi_1,\phi_2$ put
\begin{align}\label{eq:Xi-def}
    \Xi(\phi_1,\phi_2)\ &:=\ \hat\phi_1(0)\hat\phi_2(0)-\widehat{\phi_1\phi_2}(0)\ =\ \int_{-\infty}^\infty \phi_1\int_{-\infty}^\infty \phi_2-\int_{-\infty}^\infty \phi_1\phi_2,\\
    \Upsilon(\phi_1,\phi_2)\ &:=\ \hat\phi_1(0)\phi_2(0)+\hat\phi_2(0)\phi_1(0), \label{eq:Upsilon-def}\\
    \Xi_2(\phi_1,\phi_2)\ &:=\ 4\,\Xi(\phi_1,\phi_2)+\Upsilon(\phi_1,\phi_2). \label{eq:Xi2-def}
\end{align}

\begin{restatable}{theorem}{SecondUniversalityBreaks}\label{thm:D2-universality-breaks}
    Assume the Generalized Riemann Hypothesis, as in Theorem \ref{thm:D1-universality-breaks}. Let $\cF$ be an admissible family of \emph{squarefree} levels as in Definition \ref{def:profile}, and let $\phi_1,\phi_2$ be even Schwartz functions with $\hat\phi_1,\hat\phi_2$ supported in $[-\sigma,\sigma]$ for some $\sigma < 2/3$. Then
    \begin{align}\label{eq:D2-master}
        D_2(\cF_{N},\phi_1,\phi_2)\ =\ \cD_2(R;\phi_1,\phi_2)-\frac{\kappa_1(\cF)\,\Xi_2(\phi_1,\phi_2)}{\log (R)}+O\paren{\log^{-2}(R)}.
    \end{align}
    In particular
    \begin{align}\label{eq:D2-invariant}
        \lim_{N \to \infty}\log (R)\braks{D_2(\cF_{N},\phi_1,\phi_2)-\cD_2(R;\phi_1,\phi_2)}\ =\ -\kappa_1(\cF)\,\Xi_2(\phi_1,\phi_2),
    \end{align}
    which is nonzero whenever $\Xi_2(\phi_1,\phi_2) \neq 0$ and the level has a bounded prime factor, and for two such families with $N \asymp \tilde N$,
    \begin{align}\label{eq:D2-comparison}
        \lim_{N \to \infty} \log (R)\braks{D_2(\cF_{N},\phi_1,\phi_2)-D_2(\tilde\cF_{\tilde N},\phi_1,\phi_2)}\ =\ -\Xi_2(\phi_1,\phi_2)\paren{\kappa_1(\cF)-\kappa_1(\tilde\cF)}.
    \end{align}
\end{restatable}

\begin{example}\label{ex:break}
    Let $\brax{p_n}$ run through the odd primes, and compare the families with levels $\brax{N} = \brax{p_n}$,  and $\brax{\tilde N} = \brax{2p_n}$. Both sequences are admissible, and $N \asymp \tilde N$, so they have the same symmetry type and, by \eqref{eq:D1-profile-value}, the same profile up to $O\paren{\log^{-2}(R)}$. Their invariants are
    \begin{align}
        \kappa_1(\cF)\ =\ 0, \qquad
        \kappa_1(\tilde\cF)\ =\ \frac{\log (2)}{2}+\frac{\log( 2)}{3}\ =\ \frac{5\log (2)}{6}, 
    \label{eq:break-1}
\end{align}
    since $Q_0 = Q = \emptyset$ for the first and $Q_0 = Q = \brax{2}$ for the second. Hence, by \eqref{eq:D1-comparison},
    \begin{align}
        \lim_{N \to \infty}\log (R)\braks{D_1(\cF_{N},\phi)-D_1(\tilde\cF_{\tilde N},\phi)}\ =\ \frac{5\hat\phi(0)\log (2)}{3}\ \neq\ 0
    \label{eq:break-2}
\end{align}
    whenever $\hat\phi(0) \neq 0$, thus the universality of the higher order terms in the first level density breaks down in the lower order terms. Moreover, Theorem \ref{thm:D2-universality-breaks} applies to give $(5/6)\Xi_2(\phi_1,\phi_2)\log (2)$ for the corresponding two-level limit. Taking $\phi_1 = \phi_2 = \phi$ with $\hat\phi(y) = \paren{1-\abs y/\sigma}^+$ gives $\Xi_2(\phi,\phi) = 4-2\sigma/3 > 0$, so the hypothesis $\Xi_2 \neq 0$ is easy to meet.
\end{example}

\begin{remark}\label{rmk:same-symmetry}
    The families compared above all have orthogonal symmetry, and by \eqref{eq:D1-profile-value} and \eqref{eq:D2-profile-value} their densities agree in the terms predicted by random matrix theory,
    \begin{align}
        \hat\phi(0)+\frac{\phi(0)}{2}
    \label{eq:same-symmetry-1}
\end{align}
    at the first level and, at the second,
    \begin{align}
        \paren{\hat\phi_1(0)+\frac{\phi_1(0)}{2}}\paren{\hat\phi_2(0)+\frac{\phi_2(0)}{2}}+2\int_{-\infty}^\infty \abs u\hat\phi_1(u)\hat\phi_2(u)du-2\widehat{\phi_1\phi_2}(0)-\frac{\phi_1(0)\phi_2(0)}{2}.
    \label{eq:same-symmetry-2}
\end{align}
    The differences \eqref{eq:D1-comparison} and \eqref{eq:D2-comparison} are therefore invisible to the Katz--Sarnak model, which predicts only these leading terms. They mark the first place at which the arithmetic of the family asserts itself.
\end{remark}

\begin{remark}\label{rmk:oldforms}
    The mechanism behind \eqref{eq:kappa1-def} is easy to describe. A prime $q$ contributes the term $\log (q)/q$ to $\kappa_1(\cF)$ simply because $q \mid N$ removes it from the sum over primes coprime to the level. It contributes the additional $\log (q)/(q^2-1)$ when $q \Vert N$ because then $\lambda_f(q)^2 = q^{-1}$ for every newform $f$ of level $N$. Both effects are of size $\log (q)/q$, so only bounded primes survive, and a prime factor tending to infinity stays invisible at every order we compute. One plausible explanation is that the relative size of the space of oldforms is what drives the change in behavior.
\end{remark}

\begin{remark}
    The support restrictions $\sigma < 4/3$ and $\sigma < 2/3$ come only from the trivial bound \eqref{eq:normalized-general-trace-bound} for the harmonically weighted trace, and we make no attempt to optimize them. The arithmetic of the level plays no role in either one. The two-level restriction is the one-level restriction applied to $\phi_1\phi_2$, whose Fourier transform $\hat\phi_1 * \hat\phi_2$ is supported in $[-2\sigma,2\sigma]$.
\end{remark}

%%%%%%%%%%%%%%%%%%%%%%%%%%%%%%%%%%%%%%%%%%%%%%%%%%%%%%%%%%%%%%%%%%%%%%
\subsection{Structure of the argument} \label{subsec:tools}

Inclusion-exclusion writes $D_2(f,\phi_1,\phi_2)$ in terms of $D_1(f,\phi_i)$ as
\begin{equation}\label{eq:2lvl_break_down}
    D_2(f; \phi_1, \phi_2)\ =\ D_1(f, \phi_1)D_1(f, \phi_2)-2D_1(f, \phi_1\phi_2)+\frac{1-\varepsilon_f}{2}\phi_1(0)\phi_2(0),
\end{equation}
where $\varepsilon_f = \pm1$ according as the functional equation of the completed $L$-function $\Lambda(f,s)$ is even or odd.
We assume that the reader is familiar with standard properties of $L$-functions (see for example \cite{iwaniec-kwalski2004analytic}), and in particular with Satake parameters. We denote the Satake parameters of $f$ at $n$ by $\alpha_f(n)$ and $\beta_f(n)$. For $p \nmid N$ we have
\begin{equation}\label{eq:satake-good}
    \lambda_f(p) \ = \ \alpha_f(p) + \beta_f(p), \qquad
        \alpha_f(p) \beta_f(p) \ =\ 1, \qquad |\alpha_f(p)|\ =\ 1,
\end{equation}
while for $p \mid N$ we take $\alpha_f(p) \ =\ \lambda_f(p)$ and $\beta_f(p) \ =\ 0$. Throughout we abbreviate
\begin{equation}\label{eq:cm-def}
    c_m(p;f)\ :=\ \alpha_f(p)^m+\beta_f(p)^m,
\end{equation}
so that $|c_m(p;f)| \leq 2$ for all $m \geq 0$ and all $p$, and $c_m(p;f) = \lambda_f(p)^m$ when $p \mid N$.

Our main tool is the explicit formula (see \cite[(4.11)]{ILS}, for example),
\begin{equation}\label{eq:1lvl density breakdown}
    D_1(f, \phi) \ = \ \frac{U_{k,N}(\phi)}{\log(R)} - 2S(\phi;f),
    \qquad
    S(\phi;f)\ :=\ \sum_{p}\sum_{m=1}^\infty \frac{c_m(p;f)}{p^{m/2}} \frac{\log(p)}{\log(R) } \hat\phi\left(m \frac{\log(p)}{\log(R) }\right),
\end{equation}
where\footnote{ \cite{Miller2009_LOTerms_1LevelDensity} uses the notation $A_{k,N}$ instead of $U_{k,N}$.}
    \begin{equation} \label{eq:Ukn_def}
        U_{k,N} (\phi)\ :=\  2\hat\phi(0)\log\left(\frac{\sqrt{N}}{\pi} \right) + \int_{-\infty}^\infty \psi\left(\frac{k}{4}+ \frac{2\pi i x}{\log(R) }\right)\phi(x)dx + \int_{-\infty}^\infty \psi\left(\frac{k}{4}+ \frac{1}{2} +\frac{2\pi i x}{\log(R)}\right)\phi(x)dx,
    \end{equation}
and $\psi$ is the digamma function. Averaging, we set
\begin{equation} \label{eq:S1_def}
    S_1(\cF_{N},\phi)\ := \ -2\angles{S(\phi)}_N,
    \qquad
    S_2(\cF_{N},\phi_1,\phi_2)\ :=\ \angles{S(\phi_1)S(\phi_2)}_N, \end{equation}
so that, by \eqref{eq:1lvl density breakdown} and \eqref{eq:2lvl_break_down},
\begin{equation} \label{eq:D1explicitformula}
    D_1(\cF_{N},\phi) \ = \ \frac{U_{k,N}(\phi)}{\log(R)} + S_1(\cF_{N},\phi),
\end{equation}
\begin{align} \label{eq:D2explicitformula}
  D_2\left(\cF_{N}, \phi_1,\phi_2\right) \ =& \  \frac{U_{k,N}(\phi_1)U_{k,N}(\phi_2)}{\log^2(R)} +\frac{U_{k,N}(\phi_{2})}{\log(R)}S_1(\cF_{N},\phi_1)
  +\frac{U_{k,N}(\phi_{1})}{\log(R)}S_1(\cF_{N},\phi_2) \notag \\ &+4S_2(\cF_{N},\phi_1,\phi_2)- 2\left(\frac{U_{k,N}(\phi_1\phi_2)}{\log(R)} + S_1(\cF_{N},\phi_1\phi_2)\right)  \notag \\ \ &+\frac{\phi_1(0)\phi_2(0)}{2}\paren{1-\cW(\cF_{N})},
  \qquad
  \cW(\cF_{N})\ :=\ \angles{\varepsilon_f}_N.
\end{align}
The factor $4$ in front of $S_2$ arises because each $D_1$ in \eqref{eq:2lvl_break_down} contributes a factor $-2$ to the prime sum.

Everything therefore reduces to $U_{k,N}(\phi)$, $S_1(\cF_{N},\phi)$ and $S_2(\cF_{N},\phi_1,\phi_2)$. Section \ref{chp: U_k,N compute} evaluates the archimedean piece $U_{k,N}(\phi)$. Because $R = (64\pi^2)^{-1}N(k+1)(k+3)$ and $k$ is fixed, that piece depends on $N$ only through $R$, and so contributes nothing to the breaking of universality. The rest of the paper handles $S_1$ and $S_2$ as follows.

\begingroup\sloppy\begin{enumerate}[label=(\roman*)]
    \item \emph{Separation by divisibility} (Sections \ref{chp: S1 compute chpt} and \ref{chp: S2 compute chp}). Each of $S_1$ and $S_2$ splits, exactly and without error, according to whether the primes involved divide the level. At a prime $p \nmid N$ the Satake powers satisfy $c_m(p;f) = \lambda_f(p^m)-\lambda_f(p^{m-2})$, while at $p \mid N$ they collapse to $\lambda_f(p)^m$. These two identities drive everything. For $S_2$ we prove in addition that the \emph{covariance}
    \begin{align}\label{eq:cov-intro}
        \Cov(\cF_{N};\phi_1,\phi_2)\ :=\ \angles{S(\phi_1)S(\phi_2)}_N-\angles{S(\phi_1)}_N\angles{S(\phi_2)}_N
    \end{align}
    is \emph{diagonal}, meaning that up to a power-saving error only the terms with $p_1 = p_2$ survive (Theorem \ref{thm:cov-diagonal}).
    \item \emph{Uniform moment estimates} (Section \ref{chp: sum of weights and family specific terms}). We extend the Petersson--Iwaniec--Luo--Sarnak machinery to arbitrary levels and record the harmonically weighted moments of $c_m(p;f)$ and of $\lambda_f(p)^m$, uniformly in the factorization of $N$ (Lemma \ref{lem:generalized-hecke-moments} and Corollary \ref{cor:satake-moments}).
    \item \emph{Evaluation and comparison} (Section \ref{chp: lower order terms}). Substituting (ii) into (i) produces, for each of $S_1$ and $\Cov$, a \emph{universal profile}, namely the value the sum would take if no prime divided the level, together with a correction supported on the finitely many bounded primes dividing $N$. The corrections are the invariants $\kappa_1,\kappa_2,\kappa_3$ of Definition \ref{def:kappas}, and Theorems \ref{thm:D1-universality-breaks} and \ref{thm:D2-universality-breaks} follow.
\end{enumerate}

Finally, in Section \ref{sec:future-work}, we list to possible routes for future work.

%%%%%%%%%%%%%%%%%%%%%%%%%%%%%%%%%%%%%%%%%%%%%%%%%%%%%%%%%%%%%%%%%%%%%%
\section{Preliminaries}\label{chp: preliminaries}

We collect the facts about holomorphic cusp newforms that the later sections use. Nothing here is new. The two results we call on most often are the Petersson trace formula, in the forms \eqref{eq:peter-formula} and \eqref{aproxILS}, and its harmonically weighted variant \eqref{eq:harm avg}, which Section \ref{chp: sum of weights and family specific terms} extends to arbitrary level.

We define the $N^{\text{th}}$ congruence subgroup of $\SL_2(\mathbb{Z})$ by
\begin{equation}\label{eq:cong-subgroup}
\Gamma_0(N) \ := \ \left\{ \begin{pmatrix} 
a & b \\ c & d 
\end{pmatrix} \in \mathrm{SL}_2(\mathbb{Z}) \ \Bigg| \ c \equiv 0 \pmod{N} \right\}.
\end{equation}
We are interested in how the $N^{\text{th}}$ congruence subgroup acts on certain holomorphic functions $f: \mathbb{H}\to \C$, where $\mathbb{H}$ denotes the upper half complex plane. Specifically, we are interested in functions $f$ that admit a symmetry when acted on by $\Gamma_0(N)$, namely
\begin{equation}\label{eq:subsec-tools-1}
f\left( \frac{az + b}{cz + d} \right)\ = \ (cz + d)^{k} f(z).
\end{equation}
If these functions satisfy the extra stipulation that they vanish at their cusps, we call these functions $f$ a holomorphic cusp form of weight $k$ and level $N$.
Such functions admit a Fourier expansion
\begin{equation}\label{eq:Fourier}
f(z)\ =\ \sum_{n=1}^{\infty}a_f(n)e^{2\pi inz},
\end{equation}
where the coefficients $a_f(n)$ are normalized so that $a_f(1)=1$.\footnote{Since we are looking at cusp forms we have $a_f(0)=0$, which need not hold for general forms.}
We denote the space of all holomorphic cusp forms of weight $k$ and level $N$ by $S_k(N)$. This space is a finite-dimensional Hilbert space with inner product given by
\begin{equation}\label{eq:Petersson}
\langle f,g\rangle \ := \ \frac{1}{\nu(N)}\int_{\Gamma_0(N) \backslash\mathbb{H}}f(z)\overline{g(z)}\,y^k\,\frac{dx\,dy}{y^2},
\end{equation}
where $\nu(N):=[\mathrm{SL}_2(\mathbb{Z}):\Gamma_0(N)]$. This inner product is called the Petersson inner product. 
Given a modular form $f$,  we define the Hecke Operator as 
\begin{equation}\label{eq:Hecke}
T_nf(z) \ := \ n^{k-1}\sum_{\substack{ad=n \\(a,N)=1}}\sum_{b=0}^{d-1}d^{-k}f\left(\frac{az+b}{d}\right) .
\end{equation}
Importantly, for $f \in S_k(N)$,  $T_n(f(z)) \in S_k(N)$ for all $n$.
In 1970 Atkin and Lehner \cite{AL} developed the theory of newforms in $S_k(N)$. Given a form $f \in S_k(d)$ where $d\mid N$, we may also consider this to be an element of $S_k(N)$. The forms of $S_k(N)$ obtained this way are the oldforms, and they span a vector subspace $S_k^{\text{old}}(N)$. We denote its orthogonal complement under the Petersson inner product as $S_k^{\text{new}}(N)$, giving us the decomposition 
\begin{equation}\label{eq:subsec-tools-2}
    S_{k}(N) \ = \ S_k^{\text{new}}(N) \oplus S_k^{\text{old}}(N).
\end{equation}
A Hecke eigenform $f\in S_k^{\text{new}}(N)$ is a newform, and the set of equivalence classes of such a function is one dimensional \cite{MurtyMurty1997}. 
Moreover, we can find an orthogonal basis $\mathcal{B}_k(N)$ for $S_k(N)$ comprised of newforms. See \cite{ILS} for construction of this basis for square-free level and \cite{BarrettEtAl2016arXiv} for non-square-free level.
Since $S_k(N)$ is finite-dimensional, we know $|\mathcal{B}_k(N)|$ is finite \cite[(2.73)]{ILS}. Given a Hecke eigenform $f\in\mathcal{B}_k(N)$, we refer to its eigenvalue under $T_n$ as the $n^{\text{th}}$ Hecke eigenvalue of $f$, $\lambda_f(n)$. By the work of Deligne \cite{Deligne1974}, it is known that $\lambda_f(p)\in [-2,2]$ for $p \nmid N$. Moreover, the Hecke eigenvalues of $f$ are closely related to its Fourier coefficients, satisfying the relation 
\begin{equation}\label{eq:Fourier-Hecke}
    a_f(n) \ = \ \lambda_f(n)n^{(k-1)/2}.
\end{equation}
 Furthermore, the Hecke eigenvalues of $f$ possess the useful multiplicative property
\begin{align}\label{eq:mult-prop}
\lambda_f(m)\lambda_f(n)\ =\sum_{\substack{d\mid(m,n)\\ (d,N)=1}}\lambda_f\left(\frac{mn}{d^2}\right),\end{align}
so that if $(m,n)=1$, then
\begin{align}\lambda_f(mn) \ = \ \lambda_f(m)\lambda_f(n).
\label{eq:subsec-tools-3}
\end{align}
We thus define the $L$-function associated to $f$ by
\begin{equation} \label{eq:lfndefn}
L(s,f) \ := \ \sum_{n=1}^{\infty}\frac{\lambda_f(n)}{n^s}, \quad \Re(s)>1.
\end{equation}

For a cusp form $f$, we define the normalized Fourier coefficients by 
\begin{equation}\label{eq:subsec-tools-4}
    \Psi_f(n) \ := \ \left( \frac{\Gamma(k-1)}{(4\pi)^{k-1}}\right)^{1/2}\|f\|^{-1}\lambda_f(n),
\end{equation}
where $\|f\|^2:=\langle f,f\rangle$, and $\Gamma$ is the gamma function.  We wish to consider the sum
\begin{equation}\label{eq:subsec-tools-5}
    \Delta_{k,N}(m,n) \ = \ \sum_{f \in \mathcal{B}_k(N)}\bar{\Psi_f(m)}\Psi_f(n)
\end{equation}
where $\mathcal{B}_k(N)$ is an orthonormal basis of $S_k(N)$. This sum is computed using the Petersson trace formula \cite{Petersson1932, iwaniec-kwalski2004analytic}. 
\begin{prop}\label{eq:peter} We have
    \begin{equation}\label{eq:peter-formula}
        \Delta_{k,N}(m,n)\ =\ \delta(m,n)+2\pi i^k\sum_{c \equiv 0 \text{ mod}\: N}c^{-1}S(m,n;c)J_{k-1}\left( \frac{4\pi \sqrt{mn}}{c}\right)
         \end{equation}
    where $\delta(m,n)=1$ if $m \ = \ n$ and $\delta(m,n)=0$ otherwise, $J_{k-1}\left( x\right)$ denotes the Bessel function, and
    $S(m,n;c)$ is the classical Kloosterman sum.
\end{prop}
Effective bounding yields the following estimation (see \cite[Corollary 2.2]{ILS}).
\begin{prop}\label{prop:peter aprox. ILS2.2}
For $m,n \geq 1$,
        \begin{equation}\label{eq:peter-aprox-ils2-2-1}
        \Delta_{k,N}(m,n)\ =\ \delta(m,n)+O\left(\frac{\tau(N)}{k^{5/6}N}\frac{(m,n,N)\tau_3((m,n))}{((m,N)+(n,N))^{1/2}} \left(\frac{mn}{(mn)^{1/2}+kN} \right)^{1/2}\log 2mn\right),
        \end{equation}
    where the implied constant is absolute and $\tau_3(\ell):=\#\{(a,b,c)\; | \; abc=\ell \}$
\end{prop}
Moreover, we define
\begin{align}\label{def:Z_func}
Z(s, f)\ & := \ \sum_{n=1}^{\infty} \frac{\lambda_f\left(n^2\right)}{n^s}\ =\ \frac{\zeta_N(s) L(s, f \otimes f)}{\zeta(s)},
\\\label{def:Z_N_func}
Z_N(s,f)\ & :=\ \sum_{n\mid N^\infty} \frac{\lambda_f(n^2)}{n^s}.
\end{align}
We wish to evaluate the arithmetically weighted sum 
\begin{align}\label{def:Deltastar}
    \Delta^*_{k,N}(m,n)\ =\ \sum_{f \in H^*_k(N)}\frac{Z_N(1,f)}{Z(1,f)}\lambda_f(m)\lambda_f(n).
\end{align}
Iwaniec, Luo, and Sarnak \cite{ILS} derived the following expression in the case of $N$ squarefree.
\begin{prop} \label{prop:ILS_Lem_2.7}
    Let $N$ be squarefree, $(m,N)=1$, and $(n,N^2)\mid N$. Then
    \begin{equation}\label{eq:harm avg}
        \Delta^*_{k,N}(m,n)\ =\ \frac{k-1}{12}\sum_{ML=N}\frac{\mu(L)M}{\nu((n,L))}\sum_{\ell |L^{\infty}}\ell^{-1}\Delta_{k,M}(m\ell^2,n).
    \end{equation}
\end{prop}
In addition, \cite[Corollary 2.10]{ILS} provides the following estimate. 
\begin{prop}\label{eq:harm avg aprox}
 Let $N$ be squarefree, $(m,N)=1$, and $(n,N^2)\mid N$. Then,
    \begin{equation} 
        \Delta^*_{k,N}(m,n)\ =\ \frac{k-1}{12}\varphi(N)\delta(m,n)+O\left(k^{1/6}(mn)^{1/4}(n,N)^{-1/2}\tau^2(N)\tau_3((m,n))\log 2mnN \right). \label{aproxILS}
    \end{equation}
\end{prop}
  As Proposition \ref{eq:harm avg aprox} requires $N$ to be squarefree, for general $N$ we use \cite[Proposition 4.1]{BarrettEtAl2016arXiv}, which removes the squarefree restriction on $N$.
    \begin{prop} \label{prop:barrett}
    Let $(m,N)=1$ and $(n,N)=1$. Then
 \begin{equation}\label{eq:barrett-1} 
        \Delta^*_{k,N}(m,n)\ =\ \frac{k-1}{12}\sum_{LM=N}\mu(L)M\prod_{p^2\mid M}\left( \frac{p^2}{p^2-1}\right)^{-1}\sum_{\substack{\ell \mid L^{\infty}\\(\ell,M)=1}}\ell^{-1}\Delta_{k,M}(m\ell^2,n).
    \end{equation}
\end{prop}
%%%%%%%%%%%%%%%%%%%%%%%%%%%%%%%%%%%%%%%%%%%%%%%
%%%%%%%%%%%%%%%%%%%%%%%%%%%%%%%%%%%%%%%%%%%%%%

\section{Computing \texorpdfstring{${U_{k,N}(\phi)}$}{A{k,N}(phi)}} \label{chp: U_k,N compute}
We start the computation of the 1st and 2nd level density by computing ${U_{k,N}(\phi)}$ up to $O\paren{\log^{-4}(R)}$ error, which becomes the $O(\log^{-5}(R))$ error in Theorem \ref{thm:D1-universality-breaks} when divided by $\log(R)$ in \eqref{eq:D1explicitformula}. 
\begin{theorem}\label{gammaf}  Let $\phi$ be an even Schwartz function whose Fourier transform has compact support and let $\psi$ be the digamma function. Then,
    \begin{align}
U_{k,N} (\phi)& \ = \ \hat\phi(0)\log(N) +\hat\phi(0)  \left(\psi\left(\frac{k}{4}\right)+\psi\left(\frac{k}{4}+ \frac{1}{2}\right)-2\log(\pi)\right) \notag \\ &~~~~- \frac{2\pi^2}{\log^2(R)}\left(\int_{-\infty}^\infty \phi(x)x^2 dx\right)\left(\psi''\left(\frac{k}{4}\right)+\psi''\left(\frac{k}{4}+ \frac{1}{2}\right)\right) +O\paren{\log^{-4}(R)}. \label{eq:Ukn_asymp}
\end{align} 
\end{theorem}
\textit{Proof of Theorem \ref{gammaf}.} We start by proving a lemma.
\begin{lemma}  Let $\phi$ be an even Schwartz function whose Fourier transform has compact support and let $\psi$ be the digamma function. Then
    \begin{equation}\label{eq:gammaf-1}
        \int_{-\infty}^\infty \phi(x)\psi\left(\frac{k}{4}+ \frac{2\pi i x}{\log(R)}\right) dx \ =\ \psi\left(\frac{k}{4}\right)\hat\phi(0) -\frac{\psi^{''}(\frac{k}{4})2\pi^2}{\log^2(R)}\int_{-\infty}^\infty \phi(x)x^2 dx + O\paren{\log^{-4}(R)}
    \end{equation}
\end{lemma}
\begin{proof}
Using (8.363.3) of \cite{GradshteynRyzhik1965} we get that the leading term is $\psi(k/4)\hat\phi(0)$. Therefore, we subtract the leading term and compute the error term more accurately. 
Using the fact that $\phi$ is even
\begin{align}
    &\int_{-\infty}^\infty \phi(x)\psi\left(\frac{k}{4}+ \frac{2\pi i x}{\log(R)}\right) dx - \psi\left(\frac{k}{4}\right)\hat\phi(0)\notag\\ & \ = \
    \frac{1}{2}\int_{-\infty}^\infty \phi(x)\left[\psi\left(\frac{k}{4}+ \frac{2\pi i x}{\log(R)}\right) + \psi\left(\frac{k}{4}- \frac{2\pi i x}{\log(R)}\right) - 2\psi\left(\frac{k}{4}\right)\right]dx. \label{eq:goal 1 what we want di gamma expand}
\end{align}

Suppose $f:\Omega\to \C,$ where $\Omega$ is an open subset of $\C$. Suppose further that $z_0 \in \Omega$. Taylor expanding about $z_0 \in \Omega$ yields
\begin{equation}\label{eq:symm Taylor expansion up to x^4 term}
    f(z_0+z) + f(z_0-z) - 2f(z_0) \ =\  f''(z_0)z^2 + O(z^4)
\end{equation} for all $z$ such that $|z|<K$ where $K$ is the radius of convergence. Moreover, we know from Cauchy's Inequality that $K \geq \text{dist}(z_0,\del\Omega)$. 
We apply the formula in \eqref{eq:symm Taylor expansion up to x^4 term} with $\Omega = \{z\in\C : \text{Re}(z) > 0\}$, $f(z) = \psi(z)$, $z_0 = k/4$, and $z = 2\pi ix/\log(R)$. We remark that $\dist(k/4,\del(\{z\in\C : \Re(z) > 0\})) = k/4$ and hence $K\geq k/4$. Then, we get
\begin{equation}\label{eq:digamma taylor expand}
    \psi\left(\frac{k}{4} + z\right) +  \psi\left(\frac{k}{4} - z\right) - 2 \psi\left(\frac{k}{4}\right) \ =\  \psi^{''}\left(\frac{k}{4}\right) z^2 + O(z^4)
\end{equation} for all $x$ such that $|x|<k\log(R)/(8\pi)$.

Now we integrate \eqref{eq:digamma taylor expand} against the test function $\phi$ for $|x| < k\log(R)/(8\pi)$ and \eqref{eq:goal 1 what we want di gamma expand}. Then
\begin{align}\label{eq: digamma substution}
&\int_{-k\log(R)/(8\pi)}^{k\log(R)/(8\pi)}\phi(x)\psi\left(\frac{k}{4}+ \frac{2\pi i x}{\log(R)}\right) dx - \psi\left(\frac{k}{4}\right)\hat\phi(0) \notag \\ & \ = \ \frac{1}{2}\int_{-k\log(R)/(8\pi)}^{k\log(R)/(8\pi)} \phi(x)\left[-\psi''\left(\frac{k}{4}\right)\left(\frac{2\pi x}{\log(R)}\right)^2 + O\left(\frac{(2\pi x)^4}{\log^4(R)}\right)\right]dx \notag \\
& \ = \ \frac{1}{2}\int_{-k\log(R)/(8\pi)}^{k\log(R)/(8\pi)} \phi(x)\left[-\psi''\left(\frac{k}{4}\right)\left(\frac{2\pi x}{\log(R)}\right)^2\right]dx + O\paren{\log^{-4}(R)}\left(\int_{-k\log(R)/(8\pi)}^{k\log(R)/(8\pi)} \phi(x)x^4dx\right) \notag \\ 
& \ = \ \frac{1}{2}\int_{-k\log(R)/(8\pi)}^{k\log(R)/(8\pi)} \phi(x)\left[-\psi''\left(\frac{k}{4}\right)\left(\frac{2\pi x}{\log(R)}\right)^2\right]dx + O\paren{\log^{-4}(R)}.
\end{align}

We now show the tail part of the integral is negligible. Consider the series expansion of $\psi(z)$ from \cite{abramowitz1972handbook}, valid for all $z \notin \Z^{\leq 0}$,
\begin{equation}\label{eq:gammaf-2}
\psi(z) \ = \ -\gamma + \sum_{n=0}^\infty \frac{z-1}{(n+1)(n+z)}
\end{equation}
where $\gamma$ is the Euler–Mascheroni constant. For $z = k/4+2\pi i x/ \log(R)$, we then have for all $n \geq 1$ that
\begin{align}\label{ineq:helping convergence of tail}
  \left|  \frac{\left(\frac{k}{4}+2\pi i x/ \log(R)\right)-1}{(n+1)(n+\left(\frac{k}{4}+2\pi i x/ \log(R)\right))}\right| & \ = \ \left|\frac{8 i \pi x + (-4 + k) \log(R)}{(1 + n)\left(8 i \pi x + (k + 4n) \log(R)\right)}\right| \notag \\
  & \ \leq \ \frac{8\pi x + |k-4| \log(R)}{(1 + n)\left((k + 4n) \log(R)\right)} \notag \\
  & \ \leq \ \frac{8\pi x+|k-4|\log(R)}{n^2}.
\end{align} 
We also know that because $\phi(x)$ is Schwartz, $\int_{k\log(R)/(8\pi)}^\infty \phi(x) P(x)dx = O\paren{\log^{-B}(R)}$ for any polynomial $P(x)$. 
Therefore, using this fact along with \eqref{ineq:helping convergence of tail}, we find 
\begin{align}\label{eq: remainder part 1/R} &\left|\int_{k\log(R)/(8\pi)}^\infty  \phi(x)\left(-\gamma + \sum_{n=0}^\infty \frac{z-1}{(n+1)(n+z)}\right) dx \right|\notag\\
   & \ \leq \ (-\gamma +1)\left|\int_{k\log(R)/(8\pi)}^\infty  \phi(x)dx \right|+  \sum_{n=1}^\infty\left|\int_{k\log(R)/(8\pi)}^\infty\phi(x)\left(\frac{8\pi x+|k-4|\log(R)}{n^2}\right)dx \right| \notag\\
   & \ = \ \left(-\gamma +1+ \left(\sum_{n=1}^\infty \frac{|k-4|\log(R)}{n^2}\right)\right)\left|\int_{k\log(R)/(8\pi)}^\infty  \phi(x)dx\right| +  8\pi \left(\sum_{n=1}^\infty \frac{1}{n^2}\right)\left|\int_{k\log(R)/(8\pi)}^\infty  x\phi(x)dx\right|
   \notag \\
   & \ = \ O\paren{\log^{-4}(R)}.
\end{align} Similarly, we know that 
\begin{equation}\label{eq:gammaf-3}
    \left|\int_{-\infty}^{-k\log(R)/(8\pi)}  \phi(x)\left(-\gamma + \sum_{n=0}^\infty \frac{z-1}{(n+1)(n+z)}\right) dx \right| \ = \ O\paren{\log^{-4}(R)}.
\end{equation}
In addition, we know that \begin{equation} \label{eq: other piece of remainder}
    \int_{k\log(R)/(8\pi)}^{\infty} \phi(x)\left[-\psi''\left(\frac{k}{4}\right)\left(\frac{2\pi x}{\log(R)}\right)^2\right]dx + \int_{-\infty}^{-k\log(R)/(8\pi)} \phi(x)\left[-\psi''\left(\frac{k}{4}\right)\left(\frac{2\pi x}{\log(R)}\right)^2\right]dx
\end{equation} is  $O\paren{\log^{-4}(R)}$ because $\phi$ is Schwartz.
Therefore, combining \eqref{eq: digamma substution},  \eqref{eq: remainder part 1/R}, and \eqref{eq: other piece of remainder}, and noticing that \begin{align}\label{eq:x2phi-vs-phihat}
    \int_{-\infty}^{\infty} \phi(x)x^2dx \ = \ -\frac{\hat\phi''(0)}{4\pi^2}
\end{align} we have that
\begin{equation}\label{eq:gammaf-4}
\int_{-\infty}^\infty \phi(x)\psi\left(\frac{k}{4}+ \frac{2\pi i x}{\log(R)}\right) dx - \psi\left(\frac{k}{4}\right)\hat\phi(0) \  =  \ \frac{\psi''\left(\frac{k}{4}\right)}{2\log^2(R)}\hat\phi''(0) + O\paren{\log^{-4}(R)}   
\end{equation} as desired.
\end{proof}
The same argument with $z_0 = k/4+1/2$ gives the following lemma. 
\begin{lemma} Let $\phi$ be an even Schwartz function whose Fourier transform has compact support and let $\psi$ be the digamma function. Then
      \begin{equation}\label{eq:gammaf-5}
        \int_{-\infty}^\infty \phi(x)\psi\left(\frac{k}{4}+ \frac{1}{2}+ \frac{2\pi i x}{\log(R)}\right) dx \ =\ \psi\left(\frac{k}{4} + \frac{1}{2}\right)\hat\phi(0) +\frac{\psi^{''}(\frac{k}{4} + \frac{1}{2})}{2\log^2(R)}\hat\phi''(0) + O\paren{\log^{-4}(R)}. 
        \end{equation}
    \end{lemma}
Adding the two lemmas and substituting into \eqref{eq:Ukn_def} gives Theorem \ref{gammaf}. By \eqref{eq:x2phi-vs-phihat}, the $\log^{-2}(R)$ term of \eqref{eq:Ukn_asymp} may equivalently be written as $\hat\phi''(0)\paren{\psi''(k/4)+\psi''(k/4+1/2)}/(2\log^2(R))$, the form in which it is used in Lemma \ref{lem:U-universal}. \qed

\section{\texorpdfstring{Separating ${S_1(\cF_{N},\phi)}$ by divisibility}{Separating S1 by divisibility}}\label{chp: S1 compute chpt}

The whole of our analysis rests on two elementary identities for the Satake power sums $c_m(p;f)$ of \eqref{eq:cm-def}. We record them first.

\begin{lemma}\label{lem:satake-identities}
    Let $f \in H^*_k(N)$ and let $p$ be prime.
    \begin{enumerate}[label=(\alph*)]
        \item If $p \mid N$, then $c_m(p;f) = \lambda_f(p)^m$ for all $m \geq 1$.
        \item If $p \nmid N$, then $c_0(p;f) = 2$, $c_1(p;f) = \lambda_f(p)$, and
        \begin{align}\label{eq:cm-hecke}
            c_m(p;f)\ =\ \lambda_f(p^m)-\lambda_f(p^{m-2}) \qquad (m \geq 2).
        \end{align}
        \item If $p \nmid N$, then for all $a,b \geq 0$,
        \begin{align}\label{eq:cc-product}
            c_a(p;f)\,c_b(p;f)\ =\ c_{a+b}(p;f)+c_{\abs{a-b}}(p;f).
        \end{align}
        \item For all $m \geq 0$ and all $p$ we have $\abs{c_m(p;f)} \leq 2$.
    \end{enumerate}
\end{lemma}

\begin{proof}
    Part (a) is immediate from $\beta_f(p) = 0$. For (b) and (c) write $\alpha = \alpha_f(p)$, $\beta = \beta_f(p)$, so $\alpha\beta = 1$ and $\abs\alpha = \abs\beta = 1$ by \eqref{eq:satake-good}. The Hecke relation \eqref{eq:mult-prop} gives $\lambda_f(p^m) = (\alpha^{m+1}-\beta^{m+1})/(\alpha-\beta)$ for $m \geq 0$, whence
    \begin{align}
        \lambda_f(p^m)-\lambda_f(p^{m-2})\ =\ \frac{\alpha^{m+1}-\beta^{m+1}-\alpha^{m-1}+\beta^{m-1}}{\alpha-\beta}\ =\ \frac{\alpha^{m-1}\paren{\alpha^2-1}-\beta^{m-1}\paren{\beta^2-1}}{\alpha-\beta},
    \label{eq:satake-identities-1}
\end{align}
    and using $\alpha\beta = 1$ to write $\alpha^2-1 = \alpha(\alpha-\beta)$ and $\beta^2-1 = -\beta(\alpha-\beta)$ yields $\alpha^m+\beta^m$. For (c), $(\alpha^a+\beta^a)(\alpha^b+\beta^b) = \alpha^{a+b}+\beta^{a+b}+\alpha^a\beta^b+\alpha^b\beta^a$ and $\alpha^a\beta^b = \alpha^{a-b}$, $\alpha^b\beta^a = \beta^{a-b}$. Part (d) follows from $\abs\alpha = \abs\beta = 1$ when $p \nmid N$, and from $\abs{\lambda_f(p)} \leq p^{-1/2} \leq 1$ when $p \mid N$ (see \eqref{eq:level-power-reduction}).
\end{proof}

Identity \eqref{eq:cm-hecke} is the reason we work with $c_m$ rather than with the powers $\lambda_f(p)^m$. It expresses the $m$-th Satake power sum as a difference of \emph{two} Hecke eigenvalues, so averaging it against the harmonic weights costs a single application of the trace formula and loses nothing combinatorially in $m$.

Following \cite{Miller2009_LOTerms_1LevelDensity}, for a level $N$ and a prime $p$ define\footnote{\cite{Miller2009_LOTerms_1LevelDensity} denotes these sums $A_{r,\cF}$, $A'_{r,\cF}$, and denotes the two-prime sums of Section \ref{chp: S2 compute chp} similarly.}
\begin{align}
    A'_{r,N}(p)\ &:=\ \begin{cases} \angles{\lambda_f(p)^r}_N, & p \mid N, \\ 0, & p \nmid N,\end{cases} \label{def:A'_r}
    \\[4pt]
    \cA_{m,N}(p)\ &:=\ \begin{cases} \angles{c_m(p;\,\cdot\,)}_N, & p \nmid N, \\ 0, & p \mid N.\end{cases} \label{def:A_r}
\end{align}
Since the condition $p \mid N$ depends only on the level and not on the individual form, exactly one of $A'_{r,N}(p)$, $\cA_{r,N}(p)$ is a nonempty average for each $p$, and by Lemma \ref{lem:satake-identities}(a) the two agree in the sense that $\angles{c_m(p)}_N$ equals $A'_{m,N}(p)$ when $p \mid N$. This gives an exact splitting of $S_1$.

\begin{theorem}\label{theorem: S_1}
    Let $\phi$ be even and Schwartz with $\hat\phi$ supported in $[-\sigma,\sigma]$. For every level $N$,
    \begin{align}\label{eq:S1-split}
        S_1(\cF_{N},\phi)\ =\ S_{A'}(\cF_{N},\phi)+S_{A}(\cF_{N},\phi),
    \end{align}
    where
    \begin{align}
        S_{A'}(\cF_{N},\phi)\ &:=\ -2\sum_{p \mid N}\frac{\log (p)}{\log (R)}\sum_{m=1}^\infty \frac{A'_{m,N}(p)}{p^{m/2}}\hat\phi\paren{\frac{m\log (p)}{\log (R)}}, \label{eq:SA'_defn}\\
        S_{A}(\cF_{N},\phi)\ &:=\ -2\sum_{p \nmid N}\frac{\log (p)}{\log (R)}\sum_{m=1}^\infty \frac{\cA_{m,N}(p)}{p^{m/2}}\hat\phi\paren{\frac{m\log (p)}{\log (R)}}, \label{eq:SA_defn}
    \end{align}
    and both sums converge absolutely.
\end{theorem}

\begin{proof}
    By \eqref{eq:S1_def} and \eqref{eq:1lvl density breakdown},
    \begin{align}
        S_1(\cF_{N},\phi)\ =\ -2\sum_p\sum_{m=1}^\infty\frac{\angles{c_m(p)}_N}{p^{m/2}}\frac{\log (p)}{\log (R)}\hat\phi\paren{\frac{m\log (p)}{\log (R)}}.
    \label{eq:theorem-s-1-1}
\end{align}
    Absolute convergence justifies interchanging the average over $f$ with the sums over $p$ and $m$, since $\abs{c_m(p;f)} \leq 2$ gives
    \begin{align}\label{eq:S1-abs-conv}
        \sum_p\sum_{m \geq 1}\frac{\abs{c_m(p;f)}}{p^{m/2}}\frac{\log (p)}{\log (R)}\abs{\hat\phi\paren{\frac{m\log (p)}{\log (R)}}}\ \lesssim\ \frac{1}{\log (R)}\sum_p \frac{\log (p)}{p^{1/2}-1}\cdot \mathbb 1_{p \leq R^\sigma}\ \lesssim\ R^{\sigma/2},
    \end{align}
    which is finite. Splitting the sum over $p$ according to whether $p \mid N$ and applying Lemma \ref{lem:satake-identities}(a) to the first part gives \eqref{eq:S1-split}.
\end{proof}

There is no error term in \eqref{eq:S1-split}, since the decomposition is an identity. All approximation is deferred to the evaluation of $A'_{m,N}(p)$ and $\cA_{m,N}(p)$ in Section \ref{chp: sum of weights and family specific terms}, and to the passage from those values to asymptotics in Section \ref{chp: lower order terms}.

\section{\texorpdfstring{Separating ${S_2(\cF_{N},\phi_1,\phi_2)}$ and diagonalizing the covariance}{Separating S2 and diagonalizing the covariance}} \label{chp: S2 compute chp}

By \eqref{eq:S1_def} and \eqref{eq:1lvl density breakdown},
\begin{align} \label{eq:S2-raw}
    S_2(\cF_{N},\phi_1,\phi_2)\ =\ \sum_{p_1,p_2}\ \sum_{m_1,m_2 = 1}^\infty \frac{\angles{c_{m_1}(p_1)c_{m_2}(p_2)}_N}{p_1^{m_1/2}p_2^{m_2/2}}\frac{\log(p_1)\log(p_2)}{\log^2(R)}\hat\phi_1\paren{\frac{m_1\log(p_1)}{\log(R)}}\hat \phi_2 \paren{\frac{m_2\log(p_2)}{\log(R)}},
\end{align}
the rearrangement being legitimate because, exactly as in \eqref{eq:S1-abs-conv}, the double sum converges absolutely. Following the four possibilities for whether $p_1$ and $p_2$ divide $N$ we write
\begin{align}\label{eq:S2-split}
    S_2(\cF_{N},\phi_1,\phi_2)\ =\ S_{B''}(\cF_{N})+S_{B'}(\cF_{N})+S_{B}(\cF_{N}),
\end{align}
where $S_{B''}$ collects the terms with $p_1 \mid N$ and $p_2\mid N$, $S_{B'}$ those with exactly one of $p_1,p_2$ dividing $N$, and $S_B$ those with $p_1 \nmid N$ and $p_2 \nmid N$. In the notation
\begin{align}\label{not:B'' i,j}
    B''_{r_1,r_2,N}(p_1,p_2)&\ := \ \angles{\lambda_f(p_1)^{r_1}\lambda_f(p_2)^{r_2}}_N \qquad (p_1,p_2 \mid N), \\ \label{not:B' i,j}
     \cB'_{r,m,N}(p_1,p_2)&\ := \ \angles{\lambda_f(p_1)^{r}c_m(p_2;f)}_N \qquad (p_1 \mid N,\ p_2 \nmid N),\\\label{not:B i,j}
     \cB_{m_1,m_2,N}(p_1,p_2)&\ := \ \angles{c_{m_1}(p_1;f)c_{m_2}(p_2;f)}_N \qquad (p_1, p_2 \nmid N),
\end{align}
where we drop the subscript $N$ when the level is clear. Here $B''$ agrees with the quantity of the same name in Lemma \ref{lem:generalized-hecke-moments}, because $c_m(p;f) = \lambda_f(p)^m$ for $p \mid N$, while $\cB'$ and $\cB$ are the Satake analogues of the sums $B'$ and $B$ appearing there. In this notation
\begin{align}
    S_{B''}(\cF_{N}) \ = \  &\sum_{p_1,p_2 \mid N}\ \sum_{m_1,m_2=1}^{\infty}  \frac{B''_{m_1,m_2}(p_1,p_2)}{p_1^{m_1/2}p_2^{m_2/2}} \frac{\log(p_1)\log(p_2)}{\log^2(R)}\hat\phi_1 \paren{\frac{m_1 \log(p_1)}{\log(R)}}\hat \phi_2 \paren{\frac{m_2\log(p_2)}{\log(R)}}, \label{eq:SB''-def}\\
    S_{B'}(\cF_{N})   \ =\  & \sum_{(i,j) \in A}\ \sum_{\substack{p_1 \mid N\\ p_2 \nmid N}}\ \sum_{m_1,m_2=1}^{\infty}  \frac{\cB'_{m_1,m_2}(p_1,p_2)}{p_1^{m_1/2}p_2^{m_2/2}} \frac{\log(p_1)\log(p_2)}{\log^2(R)}\hat\phi_i \paren{\frac{m_1 \log(p_1)}{\log(R)}}\hat \phi_j \paren{\frac{m_2\log(p_2)}{\log(R)}}, \label{eq:SB'-def}\\
    S_{B}(\cF_{N})   \ =\  &\sum_{p_1,p_2 \nmid N}\ \sum_{m_1,m_2=1}^{\infty}  \frac{\cB_{m_1,m_2}(p_1,p_2)}{p_1^{m_1/2}p_2^{m_2/2}} \frac{\log(p_1)\log(p_2)}{\log^2(R)}\hat\phi_1 \paren{\frac{m_1 \log(p_1)}{\log(R)}}\hat \phi_2 \paren{\frac{m_2\log(p_2)}{\log(R)}}, \label{eq:SB-def}
\end{align}
where $A := \brax{(1,2),(2,1)}$ records the two ways of assigning the test functions to the two primes in the mixed case.

Written this way, $S_2$ is a genuine second moment, and the natural object is not $S_2$ itself but the covariance
\begin{align}\label{eq:cov-def}
    \Cov(\cF_{N};\phi_1,\phi_2)\ &:=\ \angles{S(\phi_1)S(\phi_2)}_N-\angles{S(\phi_1)}_N\angles{S(\phi_2)}_N \notag\\
    &\ =\ S_2(\cF_{N},\phi_1,\phi_2)-\tfrac14 S_1(\cF_{N},\phi_1)S_1(\cF_{N},\phi_2),
\end{align}
the second equality being \eqref{eq:S1_def}. Substituting \eqref{eq:cov-def} into \eqref{eq:D2explicitformula} and using \eqref{eq:D1explicitformula} gives the compact form of the two-level density that we shall use throughout,
\begin{align}\label{eq:D2-compact}
    D_2(\cF_{N},\phi_1,\phi_2)\ =\ &D_1(\cF_{N},\phi_1)D_1(\cF_{N},\phi_2)-2D_1(\cF_{N},\phi_1\phi_2) \notag\\
    &+4\Cov(\cF_{N};\phi_1,\phi_2)+\frac{\phi_1(0)\phi_2(0)}{2}\paren{1-\cW(\cF_{N})}.
\end{align}
Indeed, $D_1(\cF_N,\phi_1)D_1(\cF_N,\phi_2)$ reproduces the first three terms of \eqref{eq:D2explicitformula} together with $4\angles{S(\phi_1)}\angles{S(\phi_2)}$, and $4S_2 = 4\Cov+4\angles{S(\phi_1)}\angles{S(\phi_2)}$.

Formula \eqref{eq:D2-compact} is exact, and it isolates precisely what is new at the second level, because everything except $\Cov$ is built from one-level densities. The following theorem is the main structural result of this section. It says that the covariance is supported on the diagonal $p_1 = p_2$, so that the two-level problem collapses to a single sum over primes.

\begin{theorem}[The covariance is diagonal]\label{thm:cov-diagonal}
    Let $\phi_1,\phi_2$ be even Schwartz functions with $\hat\phi_1,\hat\phi_2$ supported in $[-\sigma,\sigma]$, and let $\eps>0$. Then
    \begin{align}\label{eq:cov-diagonal}
        \Cov(\cF_{N};\phi_1,\phi_2)\ =\ \sum_p \frac{\log^2(p)}{\log^2(R)}\sum_{m_1,m_2 = 1}^\infty \frac{\gamma_{m_1,m_2,N}(p)}{p^{(m_1+m_2)/2}}\hat\phi_1\paren{\frac{m_1\log(p)}{\log(R)}}\hat\phi_2\paren{\frac{m_2\log(p)}{\log(R)}} \notag\\
        +\ O_{k,\eps}\paren{\frac{R^{3\sigma/2}}{N^{1-\eps}}},
    \end{align}
    where $\gamma_{m_1,m_2,N}(p)$ is the local covariance
    \begin{align}\label{eq:local-cov}
        \gamma_{m_1,m_2,N}(p)\ :=\ \angles{c_{m_1}(p)c_{m_2}(p)}_N-\angles{c_{m_1}(p)}_N\angles{c_{m_2}(p)}_N,
    \end{align}
    and its values are, up to an error $O_{k,\eps}\paren{p^{(m_1+m_2)/4+\eps}N^{-1+\eps}}$,
    \begin{align}\label{eq:local-cov-values}
        \gamma_{m_1,m_2,N}(p)\ =\
        \begin{cases}
            2\delta_{m_1,m_2}-\delta_{m_1,1}\delta_{m_2,1}-\delta_{m_1,2}\delta_{m_2,2}-\delta_{\abs{m_1-m_2},2}, & p \nmid N,\\[6pt]
            p^{-(m_1+m_2)/2}, & p \Vert N,\ m_1,m_2 \text{ both odd},\\[6pt]
            0, & \text{otherwise}.
        \end{cases}
    \end{align}
\end{theorem}

\begin{proof}[Proof of Theorem \ref{thm:cov-diagonal}]
    Absolute convergence lets us subtract $\angles{S(\phi_1)}_N\angles{S(\phi_2)}_N$ from \eqref{eq:S2-raw} termwise, giving
    \begin{align}\label{eq:cov-termwise}
        \Cov(\cF_{N};\phi_1,\phi_2)\ =\ \sum_{p_1,p_2}\sum_{m_1,m_2\geq1}\frac{\gamma^{(p_1,p_2)}_{m_1,m_2,N}}{p_1^{m_1/2}p_2^{m_2/2}}\frac{\log(p_1)\log(p_2)}{\log^2(R)}\hat\phi_1\paren{\frac{m_1\log (p_1)}{\log (R)}}\hat\phi_2\paren{\frac{m_2\log (p_2)}{\log( R)}},
    \end{align}
    where
    \begin{align}
        \gamma^{(p_1,p_2)}_{m_1,m_2,N}\ :=\ \angles{c_{m_1}(p_1)c_{m_2}(p_2)}_N-\angles{c_{m_1}(p_1)}_N\angles{c_{m_2}(p_2)}_N.
    \label{eq:cov-diagonal-1}
\end{align}
    The diagonal $p_1 = p_2 = p$ contributes exactly the displayed sum in \eqref{eq:cov-diagonal}, and the values \eqref{eq:local-cov-values} are read off from Corollary \ref{cor:satake-moments}(c)--(d) by subtracting the product of the values in Corollary \ref{cor:satake-moments}(a)--(b).

    It remains to bound the off-diagonal. Fix $p_1 \neq p_2$ and abbreviate $\eps_r := r-2\floor{r/2} \in \brax{0,1}$. By Corollary \ref{cor:satake-moments} there are three cases.

    \emph{Case 1: $p_1,p_2 \nmid N$.} By Corollary \ref{cor:satake-moments}(a) and (e),
    \begin{align}
        \angles{c_{m_1}(p_1)c_{m_2}(p_2)}_N\ &=\ \delta_{m_1,2}\delta_{m_2,2}+O\paren{\frac{(p_1^{m_1}p_2^{m_2})^{1/4+\eps}}{N^{1-\eps}}}, \notag \\
        \angles{c_{m_i}(p_i)}_N\ &=\ -\delta_{m_i,2}+O\paren{\frac{p_i^{m_i(1/4+\eps)}}{N^{1-\eps}}},
    \label{eq:cov-diagonal-2}
\end{align}
    and $(-\delta_{m_1,2})(-\delta_{m_2,2}) = \delta_{m_1,2}\delta_{m_2,2}$, so the main terms cancel and 
    \begin{align}
        \gamma^{(p_1,p_2)}_{m_1,m_2,N} \ \ll_{k,\eps}\ (p_1^{m_1}p_2^{m_2})^{1/4+\eps}N^{-1+\eps}.
    \end{align}

    \emph{Case 2: $p_1,p_2 \mid N$.} If $p_i^2 \mid N$ for some $i$ then $\lambda_f(p_i) = 0$ and both terms vanish. If $p_1\Vert N$ and $p_2 \Vert N$, then by Lemma \ref{lem:generalized-hecke-moments}, when $m_1,m_2$ are both even the two terms are each equal to $p_1^{-m_1/2}p_2^{-m_2/2}$ and cancel exactly. Otherwise both are $O_{k,\eps}\paren{p_1^{-\floor{m_1/2}}p_2^{-\floor{m_2/2}}(p_1^{\eps_{m_1}}p_2^{\eps_{m_2}})^{1/4+\eps}N^{-1+\eps}}$.

    \emph{Case 3: $p_1 \Vert N$, $p_2 \nmid N$ (or the reverse).} By Corollary \ref{cor:satake-moments}(f), $\angles{\lambda_f(p_1)^{m_1}c_{m_2}(p_2)}_N = -p_1^{-m_1/2}\delta_{m_2,2}$ for $m_1$ even, up to the stated error, and $\angles{\lambda_f(p_1)^{m_1}}_N\angles{c_{m_2}(p_2)}_N = p_1^{-m_1/2}\paren{-\delta_{m_2,2}}$ for the same $m_1$, so the main terms again cancel. For $m_1$ odd both terms are $O_{k,\eps}\paren{p_1^{-\floor{m_1/2}}(p_1p_2^{m_2})^{1/4+\eps}N^{-1+\eps}}$. If $p_1^2 \mid N$ both vanish.

    In every case $\gamma^{(p_1,p_2)}_{m_1,m_2,N} \ll_{k,\eps}(p_1^{m_1}p_2^{m_2})^{1/4+\eps}N^{-1+\eps}$. Since $\hat\phi_i$ is supported in $[-\sigma,\sigma]$, the summand of \eqref{eq:cov-termwise} vanishes unless $p_i^{m_i} \leq R^\sigma$ for $i = 1,2$, so the off-diagonal contribution is
    \begin{align}
        \ll_{k,\eps}\ \frac{1}{N^{1-\eps}\log^2 (R)}\left(\sum_{\substack{p,\,m \geq 1\\ p^m \leq R^\sigma}}\frac{\log (p)}{(p^m)^{1/4-\eps}}\right)^{\!2}
        \ \ll_\eps\ \frac{1}{N^{1-\eps}\log^2 (R)}\left(\sum_{n \leq R^\sigma}\frac{\Lambda(n)}{n^{1/4-\eps}}\right)^{\!2}\ \ll_\eps\ \frac{R^{3\sigma/2}}{N^{1-\eps}},
    \label{eq:cov-diagonal-3}
\end{align}
    by the Chebyshev bound $\sum_{n \leq x}\Lambda(n)n^{-1/4+\eps} \ll_\eps x^{3/4+\eps}$. This is \eqref{eq:cov-diagonal}.
\end{proof}

\begin{remark}
    Because $R \asymp_k N$, the error term in \eqref{eq:cov-diagonal} is $O\paren{N^{3\sigma/2-1+\eps}}$, which is negligible in the sense of Remark \ref{rmk:1/N_dominates} as soon as $\sigma < 2/3$. This is the only place where the support of the test functions is constrained, as we discuss in Remark \ref{rmk:support}.
\end{remark}

\begin{remark}\label{rmk:no-truncation}
    In \eqref{eq:cov-diagonal} the sum over $m_1,m_2$ is finite for each $p$, because the support of $\hat\phi_i$ forces $m_i \leq \sigma\log (R)/\log (p)$.
\end{remark}

%%%%%%%%%%%%%%%%%%%%%%%%%%%%%%%%%%%%%%%%%%%%%%%%%%%%%%%%%%%%%%%%%%%%%%
\section{Formulas for family specific terms} \label{chp: sum of weights and family specific terms}

The goal here is a single table of harmonically weighted moments of Hecke eigenvalues and of Satake power sums, valid uniformly in the factorization of the level. We first extend the trace formula of \cite{ILS, BarrettEtAl2016arXiv} to arbitrary level and to indices divisible by primes dividing the level (Lemmas \ref{lem:modified-oldclass-sum} through \ref{lem:generalized-arbitrary-level-trace}), then record the resulting moment estimates (Lemma \ref{lem:generalized-hecke-moments} and Corollary \ref{cor:satake-moments}), and finally specialize them to a factorization profile (Corollary \ref{cor:factorization-profile-moments}).

%%%%%%%%%%%%%%%%%%%%%%%%%%%%%%%%%%%%%%%%%%%%%%%%%%%%%%%%%%%%%%%%%%%%%%
\subsection{Computing terms}

We use the harmonic weights\footnote{The harmonic weights are nearly constant across the family. By \cite{I1,HL,Miller2009_LOTerms_1LevelDensity} we have $w_R(f)= N^{\varepsilon}$ uniformly for $f$ in the family. If we allow ineffective constants, the factor $N^{\varepsilon}$ can be replaced by $\log(N)$ for sufficiently large $N$. There are other ways of normalizations. \cite{Miller2009_LOTerms_1LevelDensity} and our paper averages with 
$w_R(f) = \zeta_N(2)/Z(1,f) = \zeta(2)/L(1,\mathrm{sym}^2 f)$. Rouymi \cite{Rouymi2011TraceNonAnnulation} uses the Petersson-harmonic weight $w_R'(f)=\Gamma(k-1)(4\pi)^{1-k}\langle f,f\rangle_N^{-1}$.
Rankin-Selberg shows $\langle f,f\rangle_N \asymp_{k,N} L(1,\mathrm{sym}^2 f)$, so $w_R(f)\propto_{k,N} w_R'(f)$.
In addition, using this weight, Rouymi derived formulas for  level $N = p^{\alpha}$ for fixed prime $p$ and $\alpha \to \infty$. We also note that different weights exhibit different behaviors as seen in \cite{KnightlyReno2019} and \cite{dillon2025centered}.} to facilitate computing the asymptotics explicitly. Recall we are working with the weights 
\begin{align}
    w_R(f)\ =\ \frac{Z_N(1,f)}{Z(1,f)} \label{eq:wR_defn}
\end{align}
where $Z_N(s,f)$ and $Z(s,f)$ are defined in \eqref{def:Z_N_func} and \eqref{def:Z_func} respectively. We consider $\cF_{N} = H_k^*(N)$ the family of holomorphic cusp newforms with weight $k$ and level $N$. 

Recall that  our extended explicit formula is expressed in terms of weighted average moments of Hecke eigenvalues, namely 
\begin{equation}\label{eq:sum-of-weights-and-family-specific-1}
    \sum_{f \in H^{*}_k(N)}\lambda(p_1)^{r_1}\lambda(p_2)^{r_2},
\end{equation}
wherein restrictions are placed on $p_1,p_2$ relative to the level $N$, and, through \eqref{eq:cm-hecke}, in the corresponding averages of the Satake power sums $c_m(p;f)$. We compute these for all possible factorizations of $N$. Lemma \ref{lem:generalized-weight-and-trace-bound} through Corollary \ref{cor:factorization-profile-moments} summarize the results, and Corollary \ref{cor:satake-moments} records the form in which Section \ref{chp: lower order terms} uses them.

For an integer $r \ \geq \ 1$ define
\begin{equation}
\label{eq:nu-c-N1-def}
\nu(r) \ = \ r\prod_{p\mid r}\left(1+\frac1p\right), \qquad c(r) \ = \ \prod_{q^2\mid r}\frac{q^2}{q^2-1}, \qquad N_1 \ = \ \prod_{p\parallel N}p.
\end{equation}
We use $\xi_d(\ell)$ as the standard local coefficients in the orthonormal oldform basis of \cite{BarrettEtAl2016arXiv}.

\begin{lemma}\label{lem:modified-oldclass-sum}
Let $N = LM$, let $f\in H_k^*(M)$, and suppose that $D\mid N_1$ and $(b,N) = 1$. Define
\begin{equation}\label{eq:modified-oldclass-sum-1}
D_L \ = \ (D,L), \qquad D_M \ = \ (D,M).
\end{equation}
Then
\begin{align}\label{eq:modified-oldclass-sum}
&\sum_{d\mid L}\xi_d(1)\sum_{\ell\mid(d,D)}\xi_d(\ell)\ell^{1/2}\lambda_f\left(\frac{Db}{\ell}\right) \notag \\
&\qquad \ = \ \lambda_f(D_Mb)\prod_{p\mid D_L}\frac{\lambda_f(p)}{p+1}\prod_{\substack{p\mid L \\p\nmid M}}\rho_f(p)^{-1}\prod_{\substack{p^2\mid N \\p^2\nmid M}}\frac{p^2}{p^2-1},
\end{align}
where
\begin{equation}\label{eq:modified-oldclass-sum-2}
\rho_f(p) \ = \ 1-\frac{p\lambda_f(p)^2}{(p+1)^2}.
\end{equation}
\end{lemma}

\begin{proof}
The coefficients, $\xi_d(\ell)$, here are multiplicative in the prime-power components of both $d$ and $\ell$, so the left-hand side of \eqref{eq:modified-oldclass-sum} factors into a product of local terms. Thus at every prime which does not divide $D_L$ the local calculation is exactly the one in \cite[Lemma 3.2]{BarrettEtAl2016arXiv}. The only new calculation therefore occurs at a prime $p \mid D_L$, and there, since $D\mid N_1$, the condition $p\mid D_L$ forces $p\parallel N$, $p\parallel L$ and $p\nmid M$. In particular, the local divisor $d$ can only be 1 or $p$, having that $p$ divides $D$ exactly once. Consequently, the local contribution is 
\begin{align}\label{eq:local-D-factor-start}
T_p \ = \ \lambda_f(p)+\xi_p(1)^2\lambda_f(p)+\sqrt p\,\xi_p(1)\xi_p(p).
\end{align}
By \cite[(2.25)]{BarrettEtAl2016arXiv} and using the fact that $p\nmid M$ we have that
\begin{equation}\label{eq:modified-oldclass-sum-3}
\xi_p(p) \ = \ \rho_f(p)^{-1/2}
\end{equation}
and
\begin{equation}\label{eq:modified-oldclass-sum-4}
\xi_p(1) \ = \ -\frac{\lambda_f(p)}{\sqrt p(1+1/p)}\rho_f(p)^{-1/2} \ = \ -\frac{\sqrt p\,\lambda_f(p)}{p+1}\rho_f(p)^{-1/2}.
\end{equation}
Consequently, it follows that
\begin{equation}\label{eq:modified-oldclass-sum-5}
\xi_p(1)^2 \ = \ \frac{p\lambda_f(p)^2}{(p+1)^2\rho_f(p)} \ = \ \frac{1-\rho_f(p)}{\rho_f(p)},
\end{equation}
so
\begin{equation}\label{eq:modified-oldclass-sum-6}
1+\xi_p(1)^2 \ = \ \rho_f(p)^{-1}.
\end{equation}
At the same time
\begin{equation}\label{eq:modified-oldclass-sum-7}
\sqrt p\,\xi_p(1)\xi_p(p) \ = \ -\frac{p\lambda_f(p)}{(p+1)\rho_f(p)},
\end{equation}
and substituting these identities into \eqref{eq:local-D-factor-start} yields
\begin{align}\label{eq:modified-oldclass-sum-8}
T_p
& \ = \ \frac{\lambda_f(p)}{\rho_f(p)}-\frac{p\lambda_f(p)}{(p+1)\rho_f(p)}\notag \ = \ \frac{\lambda_f(p)}{(p+1)\rho_f(p)}.
\end{align}
Thus, the usual local factor $\rho_f(p)^{-1}$ from Barrett et al. is replaced at every prime $p\mid D_L$ by
\begin{equation}\label{eq:modified-oldclass-sum-9}
\frac{\lambda_f(p)}{p+1}\rho_f(p)^{-1}.
\end{equation}
Multiplying these modified local factors over the primes that divide $D_L$, retaining the factors from \cite{BarrettEtAl2016arXiv} at every remaining prime, and observing that the primes of $D$ contained in $M$ contribute the factor $\lambda_f(D_M)$, which comes from $\lambda_f(Db/\ell)$, leaves us with \eqref{eq:modified-oldclass-sum}.
\end{proof}

\begin{lemma}\label{lem:generalized-forward-trace}
Let $N \geq \ 2$, $D\mid N_1$, and $(ab,N) = 1$, then for every factorization $LM = N$ define
\begin{equation}\label{eq:generalized-forward-trace-1}
D_L \ = \ (D,L), \qquad D_M \ = \ (D,M).
\end{equation}
such that
\begin{align}\label{eq:generalized-forward-trace}
\Delta_{k,N}(a,Db) \ = \ \frac{12}{(k-1)N}c(N)\sum_{LM=N}\frac1{D_L}\sum_{\substack{\ell\mid L^\infty \\(\ell,M)=1}}\frac{\Delta^*_{k,M}(aD_L\ell^2,D_Mb)}{\ell}.
\end{align}
\end{lemma}

\begin{proof}
Beginning with equation (3.7) of \cite{BarrettEtAl2016arXiv} we have that since $(a, N) = 1$, the divisor in the first Fourier coefficient factor must be equal to 1, so
\begin{align}\label{eq:barrett-before-specialization}
\Delta_{k,N}(a,Db)
& \ = \ \frac{12}{(k-1)\nu(N)}\sum_{LM=N}\frac{M}{\varphi(M)}\sum_{f\in H_k^*(M)}\frac{\lambda_f(a)}{Z(1,f)} \notag \\
&\qquad\times\sum_{d\mid L}\xi_d(1)\sum_{\ell\mid(d,D)}\xi_d(\ell)\ell^{1/2}\lambda_f\left(\frac{Db}{\ell}\right).
\end{align}
Now we substitute Lemma~\ref{lem:modified-oldclass-sum} into \eqref{eq:barrett-before-specialization}. After doing this it is clear to see that in equation (2.35) of \cite{BarrettEtAl2016arXiv}, we are able to rewrite the remaining arithmetic factors to obtain
\begin{align}\label{eq:forward-with-ZN}
\Delta_{k,N}(a,Db)
\ = \ \frac{12}{(k-1)N}c(N)\sum_{LM=N}\sum_{f\in H_k^*(M)}\frac{Z_N(1,f)}{Z(1,f)}\lambda_f(a)\lambda_f(D_Mb)\prod_{p\mid D_L}\frac{\lambda_f(p)}{p+1}.
\end{align}
Continuing forward, all that remains is to rewrite the local factors in \eqref{eq:forward-with-ZN}. So, let $p\mid D_L$ then $p\nmid M$, which now has that
\begin{equation}\label{eq:generalized-forward-trace-2}
\lambda_f(p)\lambda_f(p^{2j}) \ = \ \lambda_f(p^{2j+1})+\lambda_f(p^{2j-1}),
\end{equation}
where the second term is omitted for $j \ = \ 0$. Therefore
\begin{align}
\lambda_f(p)Z_p(1,f)
 \ &= \ \sum_{j\geq0}\frac{\lambda_f(p^{2j+1})}{p^j}+\sum_{j\geq1}\frac{\lambda_f(p^{2j-1})}{p^j}
 = \ \left(1+\frac1p\right)\sum_{j\geq0}\frac{\lambda_f(p^{2j+1})}{p^j},
\label{eq:generalized-forward-trace-3}
\end{align}
giving the local identity
\begin{equation}\label{eq:odd-local-series}
\frac{\lambda_f(p)}{p+1}Z_p(1,f) \ = \ \frac1p\sum_{j\geq0}\frac{\lambda_f(p^{2j+1})}{p^j}.
\end{equation}
Further, at every prime $p\mid L$ such that $p\nmid M$ and $p\nmid D_L$, we are able to keep the usual expansion
\begin{equation}\label{eq:generalized-forward-trace-4}
Z_p(1,f) \ = \ \sum_{j\geq0}\frac{\lambda_f(p^{2j})}{p^j}.
\end{equation}
Now, since $D_L$ is squarefree and $(D_L,M) \ = \ 1$, the multiplication of these local identities gives
\begin{equation}\label{eq:local-series-product}
Z_N(1,f)\lambda_f(a)\lambda_f(D_Mb)\prod_{p\mid D_L}\frac{\lambda_f(p)}{p+1} \notag \ = \ \frac{Z_M(1,f)}{D_L}\sum_{\substack{\ell\mid L^\infty \\(\ell,M)=1}}\frac{\lambda_f(aD_L\ell^2)\lambda_f(D_Mb)}{\ell}.
\end{equation}
Thus the series in \eqref{eq:local-series-product} are absolutely convergent given that only a finite number of primes are involved. Hence substituting \eqref{eq:local-series-product} into \eqref{eq:forward-with-ZN} and collecting the sum over $f\in H_k^*(M)$ gives \eqref{eq:generalized-forward-trace}.
\end{proof}

\begin{lemma}\label{lem:generalized-arbitrary-level-trace}
Let $N \geq \ 2$, $D\mid N_1$, and $(ab,N) \ = \ 1$. For every factorization $LM \ = \ N$, define $D_L \ := \ (D,L)$ and $D_M \ := \ (D,M)$. Then
\begin{equation}\label{eq:generalized-trace-shifted}
\Delta^*_{k,N}(a,Db) \ = \ \frac{k-1}{12}\sum_{LM=N}\frac{\mu(L)M}{D_L}c(M)^{-1}\sum_{\substack{\ell\mid L^\infty  \\(\ell,M)=1}}\frac{\Delta_{k,M}(aD_L\ell^2,D_Mb)}{\ell}.
\end{equation}
Equivalently,
\begin{align}\label{eq:generalized-trace-ils-form}
\Delta^*_{k,N}(a,Db) \ = \ \frac{k-1}{12}\sum_{LM=N}\frac{\mu(L)M}{\nu(D_L)}c(M)^{-1}\sum_{\substack{\ell\mid L^\infty \\(\ell,M)=1}}\frac{\Delta_{k,M}(a\ell^2,Db)}{\ell}.
\end{align}
When $D \ = \ 1$ both formulas reduce to that of \cite{BarrettEtAl2016arXiv}'s Proposition 4.1. When $N$ is squarefree, then \eqref{eq:generalized-trace-ils-form} reduces to the corresponding formula in \cite{ILS}. Related trace formulas at non-squarefree and prime-power levels appear in \cite{I1, Rouymi2011TraceNonAnnulation}.
\end{lemma}

\begin{proof}
Let $\mathcal R$ denote the right-hand side of \eqref{eq:generalized-trace-shifted}, then apply Lemma~\ref{lem:generalized-forward-trace} at level M to each term
\begin{equation}\label{eq:generalized-arbitrary-level-trace-1}
\Delta_{k,M}(aD_L\ell^2,D_Mb).
\end{equation}
While possibly surprising, this operation is valid because
\begin{equation}\label{eq:generalized-arbitrary-level-trace-2}
D_M\mid\prod_{p\parallel M}p
\end{equation}
and
\begin{equation}\label{eq:generalized-arbitrary-level-trace-3}
(aD_L\ell^2b, M)=1.
\end{equation}
Now, for a given factorization $QW \ = \ M$ define
\begin{equation}\label{eq:generalized-arbitrary-level-trace-4}
D_Q \ = \ (D,Q), \qquad D_W \ = \ (D,W).
\end{equation}
Thus, since every prime that divides $D$ occurs exactly once in $N$, the integers $D_L$, $D_Q$, and $D_W$ must be pairwise coprime and satisfy
\begin{equation}\label{eq:generalized-arbitrary-level-trace-5}
D \ = \ D_LD_QD_W.
\end{equation}
So substituting \eqref{eq:generalized-forward-trace} into the definition of $\mathcal R$ yields
\begin{align}\label{eq:inversion-expanded}
\mathcal R
\ = \ \sum_{LQW=N}\frac{\mu(L)}{D_LD_Q}\sum_{\substack{\ell\mid L^\infty \\(\ell,QW)=1}}\frac1\ell\sum_{\substack{t\mid Q^\infty \\(t,W)=1}}\frac1t\Delta^*_{k,W}(aD_LD_Q\ell^2t^2,D_Wb).
\end{align}
Fixing $W\mid N$ and $f\in H_k^*(W)$ plus combining the local factors of the two inner series prime by prime says that if $p\mid W$ then the corresponding Euler factor is already included in $Z_W(1,f)$. However, if $p\nmid W$ and $p\nmid D$ then exactly one of the $\ell$ or $t$ sums contributes the factor
\begin{equation}\label{eq:generalized-arbitrary-level-trace-6}
Z_p(1,f) \ = \ \sum_{j\geq0}\frac{\lambda_f(p^{2j})}{p^j}.
\end{equation}
If $p\mid D$ and $p\nmid W$, then $p$ belongs to exactly one of $L$ or $ Q$, and the corresponding local factor by \eqref{eq:odd-local-series} is
\begin{equation}\label{eq:generalized-arbitrary-level-trace-7}
\frac1p\sum_{j\geq0}\frac{\lambda_f(p^{2j+1})}{p^j} \ = \ \frac{\lambda_f(p)}{p+1}Z_p(1,f).
\end{equation}
Therefore, for a fixed divisor $W$ of $N$, the entire inner contribution is
\begin{align}\label{eq:ZN-DW-definition}
\mathcal Z_{N,D}(W) \ := \ \sum_{f\in H_k^*(W)}\frac{Z_N(1,f)}{Z(1,f)}\lambda_f(a)\lambda_f(D_Wb)\prod_{\substack{p\mid D \\p\nmid W}}\frac{\lambda_f(p)}{p+1}.
\end{align}
Importantly, this expression is independent of the factorization $LQ=N/W$ and therefore we are able to write
\begin{align}
\mathcal R
& \ = \ \sum_{LQW=N}\mu(L)\mathcal Z_{N,D}(W) \ = \ \sum_{W\mid N}\mathcal Z_{N,D}(W)\sum_{L\mid N/W}\mu(L).
\label{eq:generalized-arbitrary-level-trace-8}
\end{align}
Notably, the inner sum is zero unless $W \ = \ N$, in which case it is equal to 1, but when $W \ = \ N$, we have $D_W \ = \ D$ and the product in \eqref{eq:ZN-DW-definition} is empty. Thus
\begin{equation}\label{eq:generalized-arbitrary-level-trace-9}
\mathcal R \ = \ \mathcal Z_{N,D}(N) \ = \ \Delta^*_{k,N}(a,Db),
\end{equation}
which proves \eqref{eq:generalized-trace-shifted}. Finally, it suffices to derive the equivalent form \eqref{eq:generalized-trace-ils-form}. So let $p\mid D_L$. Then since $p\nmid M$, after summing over an orthonormal basis of $S_k(M)$, equation \eqref{eq:odd-local-series} gives
\begin{align}
\frac1p\sum_{j\geq0}\frac{\Delta_{k,M}(ap^{2j+1}u,v)}{p^j} \ = \ \frac1{p+1}\sum_{j\geq0}\frac{\Delta_{k,M}(ap^{2j}u,pv)}{p^j},
\label{eq:generalized-arbitrary-level-trace-10}
\end{align}
given that all displayed factors other than the indicated powers of $p$ are to be coprime to $p$. Applying this identity independently at each prime that divides $D_L$ replaces the factor $D_L^{-1}$ and the shifted index $aD_L\ell^2$ in \eqref{eq:generalized-trace-shifted} by $\nu(D_L)^{-1}$ and the index $Db$, respectively. This gives \eqref{eq:generalized-trace-ils-form} and completes the proof.
\end{proof}

\begin{remark}
    Unless otherwise stated, a moment in the lemmas below vanishes outside the range of primes for which it is defined. Thus $A'_r(p) = 0$ for $p \nmid N$, $B''_{r_1,r_2}(p_1,p_2) = 0$ unless $p_1,p_2 \mid N$, and so on.
\end{remark}

\begin{lemma}\label{lem:generalized-weight-and-trace-bound}
Let $\cF_{N} := H_k^*(N)$ and define
\begin{equation}\label{eq:generalized-weight-and-trace-bound-1}
\mathcal A(N) = \sum_{LM=N}\mu(L)Mc(M)^{-1}.
\end{equation}
Then for fixed even $k$ and all $\eps > 0$,
\begin{equation}\label{eq:generalized-weight-asymptotic}
W_R(\cF_{N}) = \frac{k-1}{12}\mathcal A(N)+O(N^\eps).
\end{equation}
The function $\mathcal A$ is multiplicative and satisfies
\begin{equation}\label{eq:generalized-weight-and-trace-bound-2}
\mathcal A(p^e) \ = \ \begin{cases}
p-1, & e=1, \\[4pt]
p^2-p-1, & e=2, \\[4pt]
p^e\left(1-\dfrac1p\right)\left(1-\dfrac1{p^2}\right), & e\geq3.
\end{cases}
\end{equation}
In particular,
\begin{equation}\label{eq:A-lower-bound}
\mathcal A(N) \ \geq \ N\prod_{p\mid N}\left(1-\frac1p\right)^2 \gg_\eps N^{1-\eps},
\end{equation}
and consequently for all sufficiently large $N$,
\begin{equation}\label{eq:weight-lower-bound}
W_R(\cF_{N}) \ \gg_{k,\eps} \ N^{1-\eps}.
\end{equation} 
Additionally, suppose that
\begin{equation}\label{eq:generalized-weight-and-trace-bound-3}
D\mid N_1, \qquad (ab,N)=1, \qquad (a,b)=1, \qquad Dab>1.
\end{equation}
Then
\begin{equation}\label{eq:normalized-general-trace-bound}
\frac{\Delta^*_{k,N}(a,Db)}{W_R(\cF_{N})} \ \ll_{k,\eps} \ \frac{(Dab)^{1/4+\eps}}{N^{1-\eps}}.
\end{equation}
\end{lemma}

\begin{proof}
Taking $D = a = b = 1$ in \eqref{eq:generalized-trace-shifted} yields
\begin{align}\label{eq:weight-trace-expansion}
W_R(\cF_{N}) \ = \ \frac{k-1}{12}\sum_{LM=N}\mu(L)Mc(M)^{-1}\sum_{\substack{\ell\mid L^\infty \\(\ell,M)=1}}\frac{\Delta_{k,M}(\ell^2,1)}{\ell}.
\end{align}
By Proposition~\ref{prop:peter aprox. ILS2.2} and the fact that $(\ell,M) = 1$, we have
\begin{equation}\label{eq:weight-petersson-bound}
\Delta_{k,M}(\ell^2,1) \ = \ \delta_{\ell,1}+O\left(\frac{\tau(M)\ell^{1/2+\eps}}{M}\right).
\end{equation}
The diagonal term in \eqref{eq:weight-trace-expansion} is $(k-1)\mathcal A(N)/12$. For the off-diagonal contribution, choosing $\eps > 0$ sufficiently small gives
\begin{equation}\label{eq:generalized-weight-and-trace-bound-4}
\sum_{\ell\mid L^\infty}\ell^{-1/2+\eps} \ = \ \prod_{p\mid L}\left(1-p^{-1/2+\eps}\right)^{-1} \ll_\eps L^\eps.
\end{equation}
Since $\sum_{LM=N}|\mu(L)|\tau(M) \ \ll_\eps \ N^\eps$, combining these bounds shows the total off-diagonal contribution is $O(N^\eps)$, proving \eqref{eq:generalized-weight-asymptotic}. 

Because $\mathcal A$ is a Dirichlet convolution of multiplicative functions, it is multiplicative. Direct calculation on prime powers yields the formulas for $\mathcal A(p^e)$. In all cases,
\begin{equation}\label{eq:generalized-weight-and-trace-bound-5}
\frac{\mathcal A(p^e)}{p^e} \ \geq \ \left(1-\frac1p\right)^2.
\end{equation}
Multiplying over all primes dividing $N$ and using the standard estimate $\prod_{p\mid N}(1-1/p) \gg_\eps N^{-\eps}$ yields \eqref{eq:A-lower-bound}. Equations \eqref{eq:generalized-weight-asymptotic} and \eqref{eq:A-lower-bound} imply \eqref{eq:weight-lower-bound} for all sufficiently large $N$. 

For \eqref{eq:normalized-general-trace-bound}, define $x = aD_L\ell^2$ and $y = D_Mb$ in \eqref{eq:generalized-trace-shifted}. The hypotheses imply $(x,y) = 1$. If $x = y$, then $x = y = 1$, which contradicts $Dab > 1$. Thus, the diagonal term vanishes. Proposition~\ref{prop:peter aprox. ILS2.2} yields
\begin{align}\label{eq:trace-petersson-bound}
\Delta_{k,M}(aD_L\ell^2,D_Mb) \ \ll_{k,\eps} \ \frac{\tau(M)}{M}(Dab)^{1/4+\eps}\ell^{1/2+2\eps}.
\end{align}
Substituting \eqref{eq:trace-petersson-bound} into \eqref{eq:generalized-trace-shifted}, using $D_L^{-1}\leq 1$, and applying the divisor bounds gives
\begin{equation}\label{eq:generalized-weight-and-trace-bound-6}
\Delta^*_{k,N}(a,Db) \ \ll_{k,\eps} \ (Dab)^{1/4+\eps}N^\eps.
\end{equation}
Dividing by the lower bound in \eqref{eq:weight-lower-bound} and renaming $\eps$ proves \eqref{eq:normalized-general-trace-bound}.
\end{proof}

\begin{theorem}\label{thm:factorization-sensitive-weight}
Let
\begin{equation}\label{eq:factorization-sensitive-weight-1}
N_1 \ = \ \prod_{p\parallel N}p, \qquad Q(N) \ = \ \frac{N}{\operatorname{rad}(N)},
\end{equation}
and suppose $\omega(N_1) \leq t$ for some fixed $t$. Then for fixed even $k$ and every $\eps > 0$,
\begin{equation}\label{eq:factorization-sensitive-weight}
W_R(H_k^*(N)) \ = \ \frac{k-1}{12}\mathcal A(N)+O\left(N^\eps Q(N)^{-1/2}\right) \ = \ \frac{k-1}{12}\mathcal A(N)+O\left(N^\eps\sqrt{\frac{\operatorname{rad}(N)}{N}}\right).
\end{equation}
In particular, if $N$ is squarefull, then $N_1 = 1$ and $\operatorname{rad}(N)^2 \leq N$, so
\begin{equation}\label{eq:squarefull-weight-quarter}
W_R(H_k^*(N)) \ = \ \frac{k-1}{12}\mathcal A(N)+O(N^{-1/4+\eps}).
\end{equation}
\end{theorem}

\begin{proof}
In \eqref{eq:weight-trace-expansion}, the condition $\mu(L) \neq 0$ forces $L$ to be squarefree. For $\ell\mid L^\infty$ with $(\ell,M) = 1$ and $LM = N$, every prime dividing $\ell$ divides $L$ but not $M$. Thus, $\ell$ is supported only on primes dividing $N_1$. By assumption, there are at most $t$ such primes. Applying the sharper Petersson estimate from \cite[Proposition~2.2]{BarrettEtAl2016arXiv} yields
\begin{equation}\label{eq:sharper-weight-petersson}
\Delta_{k,M}(\ell^2,1) \ = \ \delta_{\ell,1}+O\left(\frac{\ell\log(2\ell^2)}{M(\ell+kM)^{1/2}}\right).
\end{equation}
Multiplying by $M/\ell$, the off-diagonal contribution for a fixed factorization $LM = N$ is bounded by
\begin{equation}\label{eq:factorization-sensitive-weight-2}
\sum_{\substack{\ell\mid L^\infty\\(\ell,M)=1}}\frac{\log(2\ell^2)}{(\ell+kM)^{1/2}}.
\end{equation}
For fixed $0 < \eta < 1/2$, we have $(\ell+kM)^{-1/2} \ll_{k,\eta} M^{-1/2+\eta}\ell^{-\eta}$. Since $\ell$ is supported on at most $t$ primes, the Euler product gives
\begin{equation}\label{eq:factorization-sensitive-weight-3}
\sum_{\substack{\ell\mid L^\infty\\(\ell,M)=1}}\frac{\log(2\ell^2)}{(\ell+kM)^{1/2}} \ \ll_{k,t,\eta} \ M^{-1/2+\eta}.
\end{equation}
Whenever $\mu(L) \neq 0$, we have $L\mid\operatorname{rad}(N)$, so $M = N/L \geq N/\operatorname{rad}(N) = Q(N)$. Since $c(M)^{-1} \leq 1$, summing over $L$, absorbing the divisor sum into $N^\eta$, and renaming $\eta$ gives
\begin{equation}\label{eq:factorization-sensitive-weight-4}
W_R(H_k^*(N))-\frac{k-1}{12}\mathcal A(N) \ \ll_{k,t,\eps} \ N^\eps Q(N)^{-1/2},
\end{equation}
proving \eqref{eq:factorization-sensitive-weight}. If $N$ is squarefull, every prime dividing $N$ occurs to at least the second power, so $\operatorname{rad}(N)^2 \leq N$. Hence $Q(N) \geq N^{1/2}$, yielding \eqref{eq:squarefull-weight-quarter}.
\end{proof}

\begin{cor}\label{cor:prime-power-weight}
Fix $e \geq 2$ and let $N = p^e$ for prime $p$. Then
\begin{equation}\label{eq:prime-power-weight-sharpening}
W_R(H_k^*(p^e)) \ = \ \frac{k-1}{12}\mathcal A(p^e)+O\left(p^{-(e-1)/2}\right).
\end{equation}
Equivalently,
\begin{equation}\label{eq:prime-power-weight-1}
W_R(H_k^*(N)) \ = \ \frac{k-1}{12}\mathcal A(N)+O\left(N^{-(e-1)/(2e)}\right).
\end{equation}
\end{cor}

\begin{proof}
For $N = p^e$, $\mu(L) \neq 0$ implies $L \in \{1, p\}$. The corresponding levels are $M = p^e$ and $M = p^{e-1}$, with the latter giving the larger error. Therefore,
\begin{equation}\label{eq:prime-power-weight-2}
W_R(H_k^*(p^e))-\frac{k-1}{12}\mathcal A(p^e) \ \ll_{k,e} \ p^{-(e-1)/2}.
\end{equation}
\end{proof}

\begin{lemma}\label{lem:generalized-hecke-moments}
Let $\cF_{N} \ := \ H_k^*(N)$ for level $N \geq 2$, and let $r, r_1,r_2 \geq 1$. Define
\begin{equation}\label{eq:generalized-hecke-moments-1}
\varepsilon_r \ = \ r-2\left\lfloor\frac r2\right\rfloor \in \{0,1\}, \qquad b_{r,r-2j} \ = \ \binom rj-\binom r{j-1}
\end{equation}
(where $\binom r{-1} \ := \ 0$), and let $C_j \ = \ \binom{2j}{j}/(j+1)$ be the $j$-th Catalan number. The following estimates hold uniformly in the prime factorization of $N$. 

If $p\mid N$, then
\begin{equation}\label{eq:generalized-hecke-moments-2}
A'_r(p) \ = \ \begin{cases}
p^{-r/2}, & p\parallel N \text{ and } r\equiv0\pmod2, \\[5pt]
O_{k,\eps}\left(\dfrac{p^{-\lfloor r/2\rfloor+1/4+\eps}}{N^{1-\eps}}\right), & p\parallel N \text{ and } r\equiv1\pmod2, \\[9pt]
0, & p^2\mid N.
\end{cases}
\end{equation}
If $p\nmid N$, then
\begin{equation}\label{eq:generalized-hecke-moments-3}
A_r(p) \ = \ \begin{cases}
C_{r/2}+O\left(\dfrac{2^r(p^r)^{1/4+\eps}}{N^{1-\eps}}\right), & r\equiv0\pmod2, \\[9pt]
O\left(\dfrac{2^r(p^r)^{1/4+\eps}}{N^{1-\eps}}\right), & r\equiv1\pmod2.
\end{cases}
\end{equation}
For $p_1 = p_2 = p$,
\begin{equation}\label{eq:generalized-hecke-moments-4}
B''_{r_1,r_2}(p,p) \ = \ A'_{r_1+r_2}(p), \qquad B_{r_1,r_2}(p,p) \ = \ A_{r_1+r_2}(p).
\end{equation}
Let $p_1\neq p_2$ with $p_1,p_2\mid N$. If $p_1^2\mid N$ or $p_2^2\mid N$, then $B''_{r_1,r_2}(p_1,p_2) = 0$. If $p_1\parallel N$ and $p_2\parallel N$, define $d = p_1^{\varepsilon_{r_1}}p_2^{\varepsilon_{r_2}}$. Then
\begin{equation}\label{eq:generalized-hecke-moments-5}
B''_{r_1,r_2}(p_1,p_2) \ = \ \begin{cases}
p_1^{-r_1/2}p_2^{-r_2/2}, & r_1\equiv r_2\equiv0\pmod2, \\[6pt]
O\left(\dfrac{p_1^{-\lfloor r_1/2\rfloor}p_2^{-\lfloor r_2/2\rfloor}d^{1/4+\eps}}{N^{1-\eps}}\right), & \text{otherwise}.
\end{cases}
\end{equation}
Let $p_1\mid N$ and $p_2\nmid N$. If $p_1^2\mid N$, then $B'_{r_1,r_2}(p_1,p_2) = 0$. If $p_1\parallel N$, then
\begin{equation}\label{eq:generalized-hecke-moments-6}
B'_{r_1,r_2}(p_1,p_2) \ = \ \begin{cases}
p_1^{-r_1/2}C_{r_2/2}+O\left(\dfrac{2^{r_2}p_1^{-\lfloor r_1/2\rfloor}(p_1^{\varepsilon_{r_1}}p_2^{r_2})^{1/4+\eps}}{N^{1-\eps}}\right), & r_1\equiv r_2\equiv0\pmod2, \\[12pt]
O\left(\dfrac{2^{r_2}p_1^{-\lfloor r_1/2\rfloor}(p_1^{\varepsilon_{r_1}}p_2^{r_2})^{1/4+\eps}}{N^{1-\eps}}\right), & \text{otherwise}.
\end{cases}
\end{equation}
Finally, for $p_1 \neq p_2$ with $p_1,p_2\nmid N$ we have
\begin{equation}\label{eq:generalized-hecke-moments-7}
B_{r_1,r_2}(p_1,p_2) \ = \ \begin{cases}
C_{r_1/2}C_{r_2/2}+O\left(\dfrac{2^{r_1+r_2}(p_1^{r_1}p_2^{r_2})^{1/4+\eps}}{N^{1-\eps}}\right), & r_1\equiv r_2\equiv0\pmod2, \\[12pt]
O\left(\dfrac{2^{r_1+r_2}(p_1^{r_1}p_2^{r_2})^{1/4+\eps}}{N^{1-\eps}}\right), & \text{otherwise}.
\end{cases}
\end{equation}
\end{lemma}

\begin{proof}
For $p\mid N$, \cite[(2.11)]{BarrettEtAl2016arXiv} yields
\begin{equation}\label{eq:generalized-hecke-moments-8}
\lambda_f(p)^2 \ = \ \begin{cases}
p^{-1}, & p\parallel N, \\
0, & p^2\mid N.
\end{cases}
\end{equation}
Thus, when $p\parallel N$, positive powers of $\lambda_f(p)$ reduce as
\begin{equation}\label{eq:level-power-reduction}
\lambda_f(p)^r \ = \ p^{-\lfloor r/2\rfloor}\lambda_f(p)^{\varepsilon_r},
\end{equation}
and vanish when $p^2\mid N$. For even $r$, averaging \eqref{eq:level-power-reduction} yields $A'_r(p) = p^{-r/2}$. For odd $r$, $A'_r(p) = p^{-\lfloor r/2\rfloor}\Delta^*_{k,N}(1,p)/W_R(\cF_{N})$, and the estimate follows from \eqref{eq:normalized-general-trace-bound} with $D = p$ and $a = b = 1$. 

For $p\nmid N$, repeated application of the Hecke relation yields
\begin{equation}\label{eq:hecke-power-expansion}
\lambda_f(p)^r \ = \ \sum_{j=0}^{\lfloor r/2\rfloor}b_{r,r-2j}\lambda_f(p^{r-2j}).
\end{equation}
This identity follows by induction. We note that $\sum_{j=0}^{\lfloor r/2\rfloor}b_{r,r-2j} = \binom r{\lfloor r/2\rfloor} \leq 2^r$ and $b_{2j,0} = C_j$. Averaging \eqref{eq:hecke-power-expansion} produces a diagonal term only for even $r$ at $j = r/2$, with coefficient $C_{r/2}$. Remaining terms are bounded using \eqref{eq:normalized-general-trace-bound} with $D = a = 1$ and $b = p^{r-2j}$. Summing the errors establishes $A_r(p)$. 

The identical prime identities for $B''$ and $B$ are immediate. Assume $p_1 \neq p_2$ and $p_1, p_2 \mid N$. If either $p_1^2\mid N$ or $p_2^2\mid N$, the corresponding eigenvalue power vanishes, giving $B''_{r_1,r_2}(p_1,p_2) = 0$. For $p_1, p_2 \parallel N$, applying \eqref{eq:level-power-reduction} yields
\begin{equation}\label{eq:generalized-hecke-moments-9}
\lambda_f(p_1)^{r_1}\lambda_f(p_2)^{r_2} \ = \ p_1^{-\lfloor r_1/2\rfloor}p_2^{-\lfloor r_2/2\rfloor}\lambda_f\left(p_1^{\varepsilon_{r_1}}p_2^{\varepsilon_{r_2}}\right).
\end{equation}
If $r_1$ and $r_2$ are both even, the remaining eigenvalue is 1, providing the main term $p_1^{-r_1/2}p_2^{-r_2/2}$. Otherwise, applying \eqref{eq:normalized-general-trace-bound} with $D = p_1^{\varepsilon_{r_1}}p_2^{\varepsilon_{r_2}}$ and $a=b=1$ yields the error term. 

For the mixed moment with $p_1\parallel N$ and $p_2\nmid N$, combining \eqref{eq:level-power-reduction} and \eqref{eq:hecke-power-expansion} gives
\begin{equation}\label{eq:generalized-hecke-moments-10}
\lambda_f(p_1)^{r_1}\lambda_f(p_2)^{r_2} \ = \ p_1^{-\lfloor r_1/2\rfloor}\sum_{j=0}^{\lfloor r_2/2\rfloor}b_{r_2,r_2-2j}\lambda_f\left(p_1^{\varepsilon_{r_1}}p_2^{r_2-2j}\right).
\end{equation}
A diagonal term appears only when both $r_1$ and $r_2$ are even, contributing $p_1^{-r_1/2}C_{r_2/2}$. Remaining terms are bounded via \eqref{eq:normalized-general-trace-bound} with $D = p_1^{\varepsilon_{r_1}}$, $a=1$, and $b=p_2^{r_2-2j}$. Using $\sum_jb_{r_2,r_2-2j} \leq 2^{r_2}$ yields the stated error for $B'_{r_1,r_2}(p_1,p_2)$. 

Finally, for $p_1 \neq p_2$ and $p_1,p_2\nmid N$, applying \eqref{eq:hecke-power-expansion} yields
\begin{align}
\lambda_f(p_1)^{r_1}\lambda_f(p_2)^{r_2}
\ = \ \sum_{j_1=0}^{\lfloor r_1/2\rfloor}\sum_{j_2=0}^{\lfloor r_2/2\rfloor}
&b_{r_1,r_1-2j_1}b_{r_2,r_2-2j_2}\lambda_f(p_1^{r_1-2j_1})\lambda_f(p_2^{r_2-2j_2}).
\label{eq:generalized-hecke-moments-11}
\end{align}
Since the trace indices are coprime, they are equal only when both are 1. This occurs when $r_1$ and $r_2$ are even, yielding the main term $C_{r_1/2}C_{r_2/2}$. Remaining terms are bounded using \eqref{eq:normalized-general-trace-bound} with $D = 1$, $a = p_1^{r_1-2j_1}$, and $b = p_2^{r_2-2j_2}$, completing the proof via $\left(\sum_{j_1}b_{r_1,r_1-2j_1}\right)\left(\sum_{j_2}b_{r_2,r_2-2j_2}\right) \leq 2^{r_1+r_2}$.
\end{proof}

\begin{cor}[Moments of the Satake power sums]\label{cor:satake-moments}
Let $\cF_{N} := H_k^*(N)$, let $\eps > 0$, and recall $c_m(p;f)$ from \eqref{eq:cm-def} and $\varepsilon_r := r-2\floor{r/2}$. The following hold uniformly in the factorization of $N$ and in $m,m_1,m_2 \geq 1$.
\begin{enumerate}[label=(\alph*)]
    \item If $p \nmid N$, then
    \begin{equation}\label{eq:cm-avg}
        \angles{c_m(p)}_N \ = \ -\delta_{m,2}+O_{k,\eps}\!\paren{\frac{p^{m(1/4+\eps)}}{N^{1-\eps}}}.
    \end{equation}
    \item If $p \mid N$, then $\angles{c_m(p)}_N = A'_{m}(p)$, evaluated in Lemma \ref{lem:generalized-hecke-moments}.
    \item If $p \nmid N$, then
    \begin{equation}\label{eq:cmcm-avg}
        \angles{c_{m_1}(p)c_{m_2}(p)}_N \ = \ 2\delta_{m_1,m_2}-\delta_{m_1+m_2,2}-\delta_{\abs{m_1-m_2},2}+O_{k,\eps}\!\paren{\frac{p^{(m_1+m_2)(1/4+\eps)}}{N^{1-\eps}}}.
    \end{equation}
    \item If $p \mid N$, then $\angles{c_{m_1}(p)c_{m_2}(p)}_N = A'_{m_1+m_2}(p)$. In particular it equals $p^{-(m_1+m_2)/2}$ when $p \Vert N$ and $m_1+m_2$ is even, it is $O_{k,\eps}\paren{p^{-\floor{(m_1+m_2)/2}+1/4+\eps}N^{-1+\eps}}$ when $p \Vert N$ and $m_1+m_2$ is odd, and it vanishes when $p^2 \mid N$.
    \item If $p_1 \neq p_2$ and $p_1,p_2 \nmid N$, then
    \begin{equation}\label{eq:cc-cross-avg}
        \angles{c_{m_1}(p_1)c_{m_2}(p_2)}_N \ = \ \delta_{m_1,2}\delta_{m_2,2}+O_{k,\eps}\!\paren{\frac{(p_1^{m_1}p_2^{m_2})^{1/4+\eps}}{N^{1-\eps}}}.
    \end{equation}
    \item If $p_1 \Vert N$ and $p_2 \nmid N$, then
    \begin{equation}\label{eq:lam-c-avg}
        \angles{\lambda_f(p_1)^{m_1}c_{m_2}(p_2)}_N \ = \ -\delta_{m_2,2}\,p_1^{-m_1/2}\mathbb 1_{2 \mid m_1}+O_{k,\eps}\!\paren{\frac{p_1^{-\floor{m_1/2}}\paren{p_1^{\varepsilon_{m_1}}p_2^{m_2}}^{1/4+\eps}}{N^{1-\eps}}},
    \end{equation}
    while if $p_1^2\mid N$ the left-hand side vanishes.
\end{enumerate}
\end{cor}

\begin{proof}
Throughout we use \eqref{eq:normalized-general-trace-bound}, which for $D \mid N_1$, $(ab,N) = 1$, $(a,b) = 1$ and $Dab > 1$ gives
\begin{equation}\label{eq:lam-avg}
    \angles{\lambda_f(Db)}_N \ = \ \frac{\Delta^*_{k,N}(1,Db)}{W_R(\cF_{N})} \ \ll_{k,\eps} \ \frac{(Db)^{1/4+\eps}}{N^{1-\eps}},
\end{equation}
together with $\angles{\lambda_f(1)}_N = 1$.

\begin{enumerate}[(a)]
    \item By Lemma \ref{lem:satake-identities}(b), $c_m(p) = \lambda_f(p^m)-\lambda_f(p^{m-2})$ for $m \geq 2$ and $c_1(p) = \lambda_f(p)$. Applying \eqref{eq:lam-avg} with $D = a = 1$ and $b = p^m$ or $b = p^{m-2}$, the only index equal to $1$ occurs for $m = 2$ in the subtracted term, contributing $-1$. Both error terms are $\ll_{k,\eps}p^{m(1/4+\eps)}N^{-1+\eps}$.
    \item Immediate from Lemma \ref{lem:satake-identities}(a).
    \item  By Lemma \ref{lem:satake-identities}(c), $c_{m_1}(p)c_{m_2}(p) = c_{m_1+m_2}(p)+c_{\abs{m_1-m_2}}(p)$. If $m_1 = m_2$ the second term is $c_0(p) = 2$ exactly, and otherwise $\abs{m_1-m_2} \geq 1$ and part (a) applies. Adding the two contributions gives \eqref{eq:cmcm-avg}, the error being $\ll_{k,\eps}p^{(m_1+m_2)(1/4+\eps)}N^{-1+\eps}$ since $\abs{m_1-m_2} \leq m_1+m_2$.
    \item  By Lemma \ref{lem:satake-identities}(a) we have $c_{m_1}(p)c_{m_2}(p) = \lambda_f(p)^{m_1+m_2}$, and Lemma \ref{lem:generalized-hecke-moments} applies.
    \item  By Lemma \ref{lem:satake-identities}(b) the product expands into at most four terms $\pm\lambda_f(p_1^{i})\lambda_f(p_2^{j})$ with $i \in \brax{m_1,m_1-2}$ and $j \in \brax{m_2,m_2-2}$, and $\lambda_f(p_1^i)\lambda_f(p_2^j) = \lambda_f(p_1^ip_2^j)$ by \eqref{eq:mult-prop}. The index $p_1^ip_2^j$ equals $1$ exactly when $i = j = 0$, which happens only for $m_1 = m_2 = 2$, and then with sign $(-1)(-1) = +1$. All other terms are bounded by \eqref{eq:lam-avg} with $a=1$ and $b = p_1^ip_2^j$.
    \item  Since $p_1 \Vert N$, \eqref{eq:level-power-reduction} gives $\lambda_f(p_1)^{m_1} = p_1^{-\floor{m_1/2}}\lambda_f(p_1)^{\varepsilon_{m_1}}$. Expanding $c_{m_2}(p_2)$ as in (e) and using $\lambda_f(p_1^{\varepsilon})\lambda_f(p_2^j) = \lambda_f(p_1^{\varepsilon}p_2^j)$, the index equals $1$ exactly when $\varepsilon_{m_1} = 0$ and $j = 0$, i.e., when $m_1$ is even and $m_2 = 2$, contributing $-p_1^{-m_1/2}$. The remaining terms are bounded by \eqref{eq:lam-avg} with $D = p_1^{\varepsilon_{m_1}}$ (which divides $N_1$) and $b = p_2^{j}$. If $p_1^2 \mid N$ then $\lambda_f(p_1) = 0$ and the average vanishes.
\end{enumerate}
\end{proof}

\begin{remark}
    Part (a) makes precise the sense in which a prime not dividing the level is generic. In the limit $N \to \infty$ the harmonic average of $c_m(p;f)$ equals $-\delta_{m,2}$ for every $m \geq 1$, which is the $m$-th Fourier coefficient of the measure $(2/\pi)\sin^2\theta\,d\theta$ governing $\lambda_f(p) = 2\cos\theta$ in the limit. Parts (b) and (d) show how sharply this fails at a prime dividing the level. There $\lambda_f(p)^2$ is not merely equidistributed but \emph{constant}, equal to $1/p$ if $p \Vert N$ and to $0$ if $p^2 \mid N$.
\end{remark}

\begin{cor}\label{cor:factorization-profile-moments}
Let $N = \prod_{i=1}^t q_i^{e_i}$ vary through a sequence of levels with fixed $t$ and fixed exponents $e_i\geq1$, and suppose
\begin{equation}\label{eq:factorization-profile-moments-1}
q_i \ \asymp \ N^{\delta_i}, \qquad \delta_i \ \geq \ 0, \qquad \sum_{i=1}^t e_i\delta_i \ = \ 1.
\end{equation}
Then the moments at primes dividing the level from Lemma~\ref{lem:generalized-hecke-moments} separate according to the individual factorization exponents. If $e_i \geq 2$, every positive moment at $q_i$ vanishes:
\begin{equation}\label{eq:profile-repeated-prime-vanishing}
A'_r(q_i) \ = \ 0,
\end{equation}
and any $B''$ or $B'$ moment containing $q_i$ as a level-dividing prime vanishes. If $e_i = 1$, then for $m \geq 1$,
\begin{equation}\label{eq:profile-even-level}
A'_{2m}(q_i) \ = \ q_i^{-m} \asymp N^{-m\delta_i},
\end{equation}
while
\begin{equation}\label{eq:profile-odd-level}
A'_{2m+1}(q_i) \ \ll_{k,\eps} \ N^{-1-\delta_i(m-1/4)+\eps}.
\end{equation}
If $i \neq j$, $e_i = e_j = 1$, and $r,s \geq 1$, then
\begin{equation}\label{eq:profile-two-level}
B''_{r,s}(q_i,q_j) \ = \
\begin{cases}
q_i^{-r/2}q_j^{-s/2}, & r \equiv s \equiv 0 \pmod2,\\[5pt]
O\left(N^{-1-\delta_i\lfloor r/2\rfloor-\delta_j\lfloor s/2\rfloor+\frac14(\delta_i\varepsilon_r+\delta_j\varepsilon_s)+\eps}\right), & \text{otherwise}.
\end{cases}
\end{equation}
Consequently, the only primes dividing the level capable of producing nondecaying even-moment main terms along such a sequence are those with $e_i = 1$ and $\delta_i = 0$. Several fixed primes may simultaneously satisfy $q_i\parallel N$, and their even-even $B''$ interactions continue to appear as products of their local factors. If $q_i\parallel N$ is fixed, $q\nmid N$, and $r,s$ are both even, then
\begin{equation}\label{eq:profile-mixed-fixed-main}
B'_{r,s}(q_i,q) \ = \ q_i^{-r/2}C_{s/2} +O\left( \frac{2^s q_i^{-r/2}q^{s(1/4+\eps)}}{N^{1-\eps}}\right).
\end{equation}
\end{cor}

\begin{proof}
The vanishing in \eqref{eq:profile-repeated-prime-vanishing} follows from $\lambda_f(q_i)^2=0$ for $q_i^2\mid N$. If $e_i = 1$, $\lambda_f(q_i)^2 = q_i^{-1}$, yielding \eqref{eq:profile-even-level}. For odd moments, Lemma~\ref{lem:generalized-hecke-moments} yields
\begin{equation}\label{eq:factorization-profile-moments-2}
A'_{2m+1}(q_i) \ \ll_{k,\eps} \ \frac{q_i^{-m+1/4+\eps}}{N^{1-\eps}}.
\end{equation}
Substituting $q_i \asymp N^{\delta_i}$ and renaming $\eps$ gives \eqref{eq:profile-odd-level}. The same substitution in the $B''$ estimate of Lemma~\ref{lem:generalized-hecke-moments} yields \eqref{eq:profile-two-level}, and \eqref{eq:profile-mixed-fixed-main} is its mixed even-even specialization.
\end{proof}

\begin{remark}\label{rmk:bounded-primes-only}
    Corollary \ref{cor:factorization-profile-moments} demonstrates exactly how bounded prime factors create the only non-negligible contributions. If the level has no bounded prime factor, then $A'$, $B''$ and $B'$ have no main term at all, since every one of their contributions decays like a fixed negative power of $N$. A factor $q_i^{-r/2}$ decays like $N^{-\delta_i r/2}$ with $\delta_i > 0$. On the other hand, a bounded prime factor produces main terms that do not decay, such as $A'_{2m}(q) = q^{-m}$ exactly, and these survive into the lower order terms of the densities.
\end{remark}

%%%%%%%%%%%%%%%%%%%%%%%%%%%%%%%%%%%%%%%%%%%%%%%%%%%%%%%%%
%%%%%%%%%%%%%%%%%%%%%%%%%%%%%%%%%%%%%%%%%%%%%%%%%%%%%%%%%%
%%%%%%%%%%%%%%%%%%%%%%%%%%%%%%%%%%%%%%%%%%%%%%%%%%%%%%%%%

\section{Lower order terms and the breaking of universality} \label{chp: lower order terms}

We may now substitute the moment estimates of Section \ref{chp: sum of weights and family specific terms} into the decompositions of Sections \ref{chp: S1 compute chpt} and \ref{chp: S2 compute chp}. Each of $S_1$ and $\Cov$ splits as a universal profile, defined in Subsection \ref{subsec:profiles} and evaluated in Subsection \ref{subsec:profile-values}, plus a correction supported on the bounded primes dividing the level. Subsection \ref{subsec:onelevel} computes the one-level correction and Subsection \ref{subsec:cov} the two-level one, and Subsection \ref{subsec:provingBU} assembles them into Theorems \ref{thm:D1-universality-breaks} and \ref{thm:D2-universality-breaks}.

Throughout, $\cF$ is an admissible family in the sense of Definition \ref{def:profile}, with invariants $Q_0 = Q_0(\cF)$, $Q = Q(\cF)$ and $\kappa_1,\kappa_2,\kappa_3$ as in Definition \ref{def:kappas}. We write $\delta := \min\brax{\delta_i \ : \ \delta_i > 0} > 0$ for the smallest positive exponent in \eqref{eq:level-shape}, so that every prime factor of $N$ outside $Q_0$ is $\gg N^{\delta}$.

\begin{remark} \label{rmk:1/N_dominates}
    Since $R \asymp_k N$ we have $\log (N) \sim \log (R)$, so for every $\omega > 0$ and every $C > 0$,
    \begin{align}
        \frac{1}{N^\omega} \ = \ O_{\omega,C}\paren{\log^{-C}(R)}.
    \label{eq:1-n-dominates-1}
\end{align}
    Errors that are a negative power of $N$ occur throughout, and we call such terms {negligible}. We absorb them without comment into any of the error terms $O\paren{\log^{-C}(R)}$ below.
\end{remark}

\begin{remark}\label{rmk:support}
    Two error terms in this section are power-saving in $N$ but grow with the support of the test functions. The trace error in Theorem \ref{thm:SA} is $O\paren{R^{3\sigma/4}N^{-1+\eps}}$, and the off-diagonal error of Theorem \ref{thm:cov-diagonal} is $O\paren{R^{3\sigma/2}N^{-1+\eps}}$. As $R \asymp_k N$, these are negligible precisely when $\sigma < 4/3$ and $\sigma < 2/3$ respectively. Every other error below is either $O\paren{\log^{-C}(R)}$ with the constant independent of $\sigma$, or negligible for all $\sigma$. This is the sole origin of the support restrictions in Theorems \ref{thm:D1-universality-breaks} and \ref{thm:D2-universality-breaks}.
\end{remark}

%%%%%%%%%%%%%%%%%%%%%%%%%%%%%%%%%%%%%%%%%%%%%%%%%%%%%%%%%%%%%%%%%%%%%%
\subsection{Universal profiles}\label{subsec:profiles}

We first record that the archimedean term carries no information about the factorization.

\begin{lemma}\label{lem:U-universal}
    Let $\phi$ be even and Schwartz with $\hat\phi$ compactly supported. Then $U_{k,N}(\phi)$ depends on $N$ only through $R$, and writing $U_k(R;\phi) := U_{k,N}(\phi)$ we have
    \begin{align}\label{eq:U-over-logR}
        \frac{U_{k}(R;\phi)}{\log (R)}\ =\ \hat\phi(0)+\frac{\gamma_{k,1}\hat\phi(0)}{\log (R)}+\frac{\gamma_{k,2}\hat\phi''(0)}{2\log^3 (R)}+O\paren{\log^{-5}(R)},
    \end{align}
    where
    \begin{align}\label{eq:gamma-k-def}
        \gamma_{k,1}\ &:=\ \psi\paren{\frac{k}{4}}+\psi\paren{\frac{k}{4}+\frac12}-2\log(\pi)-\log\paren{\frac{(k+1)(k+3)}{64\pi^2}}, \notag\\
        \gamma_{k,2}\ &:=\ \psi''\paren{\frac{k}{4}}+\psi''\paren{\frac{k}{4}+\frac12}.
    \end{align}
\end{lemma}

\begin{proof}
    Since $R = (64\pi^2)^{-1}N(k+1)(k+3)$ with $k$ fixed,
    \begin{align}\label{eq:logN-logR}
        \log(N)\ =\ \log(R)-\log\paren{\frac{(k+1)(k+3)}{64\pi^2}},
    \end{align}
    so every occurrence of $N$ in \eqref{eq:Ukn_def} is a function of $R$. Theorem \ref{gammaf} gives
    \begin{align}\label{eq:u-universal-1}
        U_{k,N}(\phi)\ =\ \hat\phi(0)\log( N)+\hat\phi(0)\paren{\psi\paren{\tfrac k4}+\psi\paren{\tfrac k4+\tfrac12}-2\log (\pi)}\notag\\ -\frac{2\pi^2\gamma_{k,2}}{\log^2 (R)}\int_{-\infty}^\infty x^2\phi(x)dx+O\paren{\log^{-4}(R)},
    \end{align}
    and $\int_{-\infty}^\infty x^2\phi(x)dx = -\hat\phi''(0)/(4\pi^2)$. Substituting \eqref{eq:logN-logR} and dividing by $\log (R)$ gives \eqref{eq:U-over-logR}.
\end{proof}

The remaining ingredients of $D_1$ and $D_2$ are prime sums, and each has a canonical \emph{universal profile}, namely the value the sum would take if no prime divided the level.

\begin{definition}[Universal profiles]\label{def:profiles}
    For $R > 1$ and even Schwartz $\phi,\phi_1,\phi_2$ with compactly supported Fourier transforms, set
    \begin{align}
        \cE_1(R;\phi)\ &:=\ \frac{2}{\log (R)}\sum_{p}\frac{\log( p)}{p}\hat\phi\paren{\frac{2\log (p)}{\log (R)}}, \label{eq:E1-def}\\
        \cE_2(R;\phi_1,\phi_2)\ &:=\ \sum_p\frac{\log^2 (p)}{\log^2 (R)}\sum_{m_1,m_2 = 1}^\infty\frac{\gamma^{\gen}_{m_1,m_2}}{p^{(m_1+m_2)/2}}\hat\phi_1\paren{\frac{m_1\log (p)}{\log (R)}}\hat\phi_2\paren{\frac{m_2\log (p)}{\log (R)}}, \label{eq:E2-def}
    \end{align}
    where
    \begin{align}\label{eq:gamma-gen}
        \gamma^{\gen}_{m_1,m_2}\ :=\ 2\delta_{m_1,m_2}-\delta_{m_1,1}\delta_{m_2,1}-\delta_{m_1,2}\delta_{m_2,2}-\delta_{\abs{m_1-m_2},2}
    \end{align}
    is the generic local covariance of \eqref{eq:local-cov-values}. Further set
    \begin{align}
        \cD_1(R;\phi)\ &:=\ \frac{U_k(R;\phi)}{\log (R)}+\cE_1(R;\phi), \label{eq:D1-profile-def}\\
        \cD_2(R;\phi_1,\phi_2)\ &:=\ \cD_1(R;\phi_1)\cD_1(R;\phi_2)-2\cD_1(R;\phi_1\phi_2)+4\cE_2(R;\phi_1,\phi_2)+\frac{\phi_1(0)\phi_2(0)}{2}. \label{eq:D2-profile-def}
    \end{align}
\end{definition}

Every quantity in Definition \ref{def:profiles} is a function of $k$, the test functions, and $R$ alone, and none of them depend on the factorization of $N$. This is what makes the comparisons \eqref{eq:D1-invariant} and \eqref{eq:D2-invariant} meaningful. Theorem \ref{thm:profile-values} evaluates the profiles, and their leading terms turn out to be exactly the Katz--Sarnak predictions for orthogonal symmetry.

%%%%%%%%%%%%%%%%%%%%%%%%%%%%%%%%%%%%%%%%%%%%%%%%%%%%%%%%%%%%%%%%%%%%%%
\subsection{The one-level density}\label{subsec:onelevel}

\begin{theorem} \label{thm:SA' universality}
    Let $\phi$ be an even Schwartz function with $\hat\phi$ compactly supported, and set
    \begin{align}\label{eq:beta-def}
        \beta_{1}(\cF)\ :=\ \sum_{q \in Q}\frac{\log (q)}{q^2-1},
        \qquad
        \beta_{3}(\cF)\ :=\ \sum_{q \in Q}\frac{q^2\paren{q^2+1}\log^3 (q)}{\paren{q^2-1}^3}.
    \end{align}
    Then
    \begin{align}\label{eq:SA'-explicit}
        S_{A'}(\cF_{N},\phi)\ =\ -\frac{2\hat\phi(0)\beta_{1}(\cF)}{\log(R)}-\frac{4\hat\phi''(0)\beta_{3}(\cF)}{\log^3(R)}+O\paren{\log^{-5}(R)}.
    \end{align}
\end{theorem}

\begin{proof}
    By \eqref{eq:SA'_defn} and Lemma \ref{lem:generalized-hecke-moments}, $A'_{m,N}(p) = 0$ whenever $p^2 \mid N$, so only the primes $p \Vert N$ contribute. For such $p$ the series converges absolutely, since $\abs{A'_{m,N}(p)} \leq 1$ gives $\sum_{m \geq 1}\abs{A'_{m,N}(p)}p^{-m/2} \leq (p^{1/2}-1)^{-1}$, so we may split it into even and odd $m$.

    For $m$ odd, Lemma \ref{lem:generalized-hecke-moments} gives $A'_{m,N}(p) \ll_{k,\eps} p^{-\floor{m/2}+1/4+\eps}N^{-1+\eps}$, so the odd terms contribute
    \begin{align}
        \ \ll_{k,\eps}\ \frac{1}{N^{1-\eps}}\sum_{p \Vert N}\frac{\log (p)}{\log (R)}\sum_{\substack{m \geq 1 \\ m \text{ odd}}}p^{-\floor{m/2}-m/2+1/4+\eps}
        \ \ll\ \frac{n\log (N)}{N^{1-\eps}\log (R)},
    \label{eq:sa-universality-1}
\end{align}
    which is negligible. For $m = 2\ell$ even, Lemma \ref{lem:generalized-hecke-moments} gives $A'_{2\ell,N}(p) = p^{-\ell}$ exactly, so
    \begin{align}\label{eq:SA'-even-part}
        S_{A'}(\cF_{N},\phi)\ =\ -2\sum_{p \Vert N}\frac{\log (p)}{\log (R)}\sum_{\ell = 1}^\infty \frac{1}{p^{2\ell}}\hat\phi\paren{\frac{2\ell\log (p)}{\log (R)}}+O\paren{N^{-1/2}}.
    \end{align}
    If $p \Vert N$ and $p \notin Q_0$ then $p \gg N^{\delta}$, and bounding $\hat\phi$ and $\log (p)/\log (R)$ trivially shows that the corresponding summand is $\ll p^{-2} \ll N^{-2\delta}$, again negligible. Only the primes $q \in Q$ survive.

    Fix $q \in Q$. As $\hat\phi$ is even and Schwartz, $\hat\phi'(0) = 0$ and Taylor's theorem gives, for all $y \in \bR$,
    \begin{align}\label{eq:phihat-taylor-4}
        \hat\phi(y)\ =\ \hat\phi(0)+\frac{\hat\phi''(0)}{2}y^2+E_\phi(y), \qquad \abs{E_\phi(y)}\ \leq\ \frac{\norm{\hat\phi^{(4)}}_\infty}{24}y^4 .
    \end{align}
    Setting $y = 2\ell\log (q)/\log (R)$, so that the quadratic term carries a factor $4\ell^2$, and using
    \begin{align}
        \sum_{\ell=1}^\infty \frac{1}{q^{2\ell}}\ =\ \frac{1}{q^2-1},\qquad
        \sum_{\ell=1}^\infty \frac{\ell^2}{q^{2\ell}}\ =\ \frac{q^2\paren{q^2+1}}{\paren{q^2-1}^3}, \qquad
        \sum_{\ell=1}^\infty \frac{\ell^4}{q^{2\ell}}\ <\ \infty,
    \label{eq:sa-universality-2}
\end{align}
    we obtain
    \begin{align}
        \sum_{\ell=1}^\infty \frac{1}{q^{2\ell}}\hat\phi\paren{\frac{2\ell\log (q)}{\log (R)}}
        \ =\ \frac{\hat\phi(0)}{q^2-1}+\frac{2\hat\phi''(0)\log^2 (q)}{\log^2 (R)}\cdot\frac{q^2\paren{q^2+1}}{\paren{q^2-1}^3}+O\paren{\frac{\log^4 (q)}{\log^4 (R)}}.
    \label{eq:sa-universality-3}
\end{align}
    Multiplying by $-2\log (q)/\log (R)$ and summing over the finitely many $q \in Q$ yields \eqref{eq:SA'-explicit}.
\end{proof}

\begin{theorem}\label{thm:SA}
    Let $\phi$ be an even Schwartz function with $\hat\phi$ supported in $[-\sigma,\sigma]$ for some $\sigma < 4/3$, and set
    \begin{align}\label{eq:alpha-def}
        \alpha_{1}(\cF)\ :=\ \sum_{q \in Q_0}\frac{\log (q)}{q},
        \qquad
        \alpha_{3}(\cF)\ :=\ \sum_{q \in Q_0}\frac{\log^3 (q)}{q}.
    \end{align}
    Then
    \begin{align}\label{eq:SA_asymp}
        S_{A}(\cF_{N},\phi)\ =\ \cE_1(R;\phi)-\frac{2\hat\phi(0)\alpha_{1}(\cF)}{\log (R)}-\frac{4\hat\phi''(0)\alpha_{3}(\cF)}{\log^3 (R)}+O\paren{\log^{-5}(R)}.
    \end{align}
\end{theorem}

\begin{proof}
    By Corollary \ref{cor:satake-moments}(a), for $p \nmid N$ we have $\cA_{m,N}(p) = -\delta_{m,2}+O_{k,\eps}\paren{p^{m(1/4+\eps)}N^{-1+\eps}}$. Inserting this into \eqref{eq:SA_defn}, the main term is
    \begin{align}
        -2\sum_{p \nmid N}\frac{\log (p)}{\log (R)}\cdot\frac{-1}{p}\hat\phi\paren{\frac{2\log (p)}{\log (R)}}\ =\ \frac{2}{\log (R)}\sum_{p \nmid N}\frac{\log (p)}{p}\hat\phi\paren{\frac{2\log (p)}{\log (R)}}.
    \label{eq:sa-1}
\end{align}
    For the error, note that since $\hat\phi$ is supported in $[-\sigma,\sigma]$ the summand of \eqref{eq:SA_defn} vanishes unless $p^m \leq R^\sigma$. Hence the error contributes
    \begin{align}\label{eq:SA-error}
        \ \ll_{k,\eps}\ \frac{1}{N^{1-\eps}\log (R)}\sum_{\substack{p,\,m \geq 1\\p^m \leq R^\sigma}}\frac{\log (p)}{(p^m)^{1/4-\eps}}
        \ \ll_\eps\ \frac{1}{N^{1-\eps}\log( R)}\sum_{n \leq R^\sigma}\frac{\Lambda(n)}{n^{1/4-\eps}}
        \ \ll_\eps\ \frac{R^{3\sigma/4}}{N^{1-\eps}\log (R)},
    \end{align}
    which is negligible for $\sigma < 4/3$ because $R \asymp_k N$. No truncation of the $m$-sum was required, since the compact support of $\hat\phi$ supplies one and $\cA_{m,N}(p)$ carries no factor growing with $m$ (Remark \ref{rmk:no-truncation}).

    It remains to compare $\sum_{p \nmid N}$ with $\sum_p$. By \eqref{eq:E1-def},
    \begin{align}\label{eq:SA-deletion}
        \frac{2}{\log (R)}\sum_{p \nmid N}\frac{\log (p)}{p}\hat\phi\paren{\frac{2\log (p)}{\log (R)}}\ =\ \cE_1(R;\phi)-\frac{2}{\log (R)}\sum_{p \mid N}\frac{\log (p)}{p}\hat\phi\paren{\frac{2\log (p)}{\log (R)}}.
    \end{align}
    In the subtracted sum, a prime $p \mid N$ with $p \notin Q_0$ satisfies $p \gg N^\delta$, so its summand is $\ll N^{-\delta}\log( N)/\log (R)$. There are at most $n$ such primes, so their total contribution is negligible. For the finitely many $q \in Q_0$, \eqref{eq:phihat-taylor-4} with $y = 2\log (q)/\log (R)$ gives
    \begin{align}
        \hat\phi\paren{\frac{2\log (q)}{\log (R)}}\ =\ \hat\phi(0)+\frac{2\hat\phi''(0)\log^2 (q)}{\log^2 (R)}+O\paren{\log^{-4}(R)},
    \label{eq:sa-2}
\end{align}
    so that
    \begin{align}
        \frac{2}{\log (R)}\sum_{p \mid N}\frac{\log (p)}{p}\hat\phi\paren{\frac{2\log (p)}{\log (R)}}\ =\ \frac{2\hat\phi(0)\alpha_1(\cF)}{\log (R)}+\frac{4\hat\phi''(0)\alpha_3(\cF)}{\log^3 (R)}+O\paren{\log^{-5}(R)}.
    \label{eq:sa-3}
\end{align}
    Combining the last three displays gives \eqref{eq:SA_asymp}.
\end{proof}

Adding Theorems \ref{thm:SA' universality} and \ref{thm:SA} to the exact splitting \eqref{eq:S1-split}, and observing that $\kappa_1 = \alpha_1+\beta_1$ and $\kappa_3 = \alpha_3+\beta_3$ by \eqref{eq:kappa1-def}, \eqref{eq:kappa3-def}, \eqref{eq:beta-def} and \eqref{eq:alpha-def}, gives the following.

\begin{cor}\label{cor:S1-master}
    Let $\phi$ be an even Schwartz function with $\hat\phi$ supported in $[-\sigma,\sigma]$ for $\sigma < 4/3$. Then
    \begin{align}\label{eq:S1-master}
        S_1(\cF_{N},\phi)\ =\ \cE_1(R;\phi)-\frac{2\hat\phi(0)\kappa_1(\cF)}{\log (R)}-\frac{4\hat\phi''(0)\kappa_3(\cF)}{\log^3 (R)}+O\paren{\log^{-5}(R)}.
    \end{align}
\end{cor}

% \begin{remark}
%     Corollary \ref{cor:S1-master} separates the two mechanisms of Remark \ref{rmk:oldforms}. The term $\alpha_1(\cF) = \sum_{q \in Q_0}\log (q)/q$ of Theorem \ref{thm:SA} is a \emph{deletion}. A bounded prime dividing the level is simply absent from the generic sum $\cE_1$, so its generic contribution $-2\log (q)/(q\log (R))$ must be removed. The term $\beta_1(\cF) = \sum_{q \in Q}\log (q)/(q^2-1)$ of Theorem \ref{thm:SA' universality} is a \emph{replacement}. The prime is put back, but with the rigid local behavior $\lambda_f(q)^2 = q^{-1}$ forced on newforms of level exactly divisible by $q$. A prime with $q^2 \mid N$ is deleted and never replaced, because then $\lambda_f(q) = 0$, and that is why $\kappa_1(\cF)$ detects every bounded prime factor of the level rather than only those dividing it exactly once.
% \end{remark}

%%%%%%%%%%%%%%%%%%%%%%%%%%%%%%%%%%%%%%%%%%%%%%%%%%%%%%%%%%%%%%%%%%%%%%
\subsection{The covariance}\label{subsec:cov}

\begin{theorem}\label{thm:cov-value}
    Let $\phi_1,\phi_2$ be even Schwartz functions with $\hat\phi_1,\hat\phi_2$ supported in $[-\sigma,\sigma]$ for some $\sigma < 2/3$. Then
    \begin{align}\label{eq:cov-value}
        \Cov(\cF_{N};\phi_1,\phi_2)\ =\ \cE_2(R;\phi_1,\phi_2)-\frac{\hat\phi_1(0)\hat\phi_2(0)\,\kappa_2(\cF)}{\log^2 (R)}+O\paren{\log^{-4}(R)}.
    \end{align}
\end{theorem}

\begin{proof}
    By Theorem \ref{thm:cov-diagonal} and Remark \ref{rmk:support}, up to a negligible error
    \begin{align}\label{eq:cov-two-blocks}
        &\Cov(\cF_{N};\phi_1,\phi_2)\ =\ \sum_{p \nmid N}\Gamma_p+\sum_{p \mid N}\Gamma_p, \notag \\
        &\Gamma_p\ :=\ \frac{\log^2 (p)}{\log^2 (R)}\sum_{m_1,m_2 \geq 1}\frac{\gamma_{m_1,m_2,N}(p)}{p^{(m_1+m_2)/2}}\hat\phi_1\paren{\tfrac{m_1\log( p)}{\log (R)}}\hat\phi_2\paren{\tfrac{m_2\log (p)}{\log (R)}}.
    \end{align}

    \emph{The primes $p \nmid N$.} Here $\gamma_{m_1,m_2,N}(p) = \gamma^{\gen}_{m_1,m_2}+O_{k,\eps}\paren{p^{(m_1+m_2)(1/4+\eps)}N^{-1+\eps}}$ by \eqref{eq:local-cov-values}. As in \eqref{eq:SA-error}, the compact support restricts to $p^{m_i} \leq R^\sigma$, and the total error is
    \begin{align}
        \ \ll_{k,\eps}\ \frac{1}{N^{1-\eps}\log^2 (R)}\sum_{p \leq R^\sigma}\log^2 (p)\left(\sum_{m \geq 1}p^{-m/4+\eps m}\right)^2\ \ll_\eps\ \frac{R^{\sigma}}{N^{1-\eps}},
    \label{eq:cov-value-1}
\end{align}
    which is negligible. Hence, writing $\Gamma^{\gen}_p$ for $\Gamma_p$ with $\gamma_{m_1,m_2,N}(p)$ replaced by $\gamma^{\gen}_{m_1,m_2}$,
    \begin{align}\label{eq:cov-generic-part}
        \sum_{p \nmid N}\Gamma_p\ =\ \cE_2(R;\phi_1,\phi_2)-\sum_{p \mid N}\Gamma^{\gen}_p+O\paren{N^{-\eta}}
    \end{align}
    for some $\eta > 0$, by \eqref{eq:E2-def}.

    \emph{Evaluating $\Gamma^{\gen}_q$ for $q$ bounded.} Every pair $(m_1,m_2) \neq (1,1)$ with $\gamma^{\gen}_{m_1,m_2} \neq 0$ has $m_1+m_2 \geq 4$, so replacing each $\hat\phi_i$ by $\hat\phi_i(0)$ costs $O\paren{\log^{-4}(R)}$ in total by \eqref{eq:phihat-taylor-4} and the convergence of $\sum_{m_1,m_2}(m_1^2+m_2^2)q^{-(m_1+m_2)/2}$. Since
    \begin{align}\label{eq:Gq}
        G(q)\ :=\ \sum_{m_1,m_2 \geq 1}\frac{\gamma^{\gen}_{m_1,m_2}}{q^{(m_1+m_2)/2}}
        &\ =\ \underbrace{\sum_{m \geq 1}\frac{2-\delta_{m,1}-\delta_{m,2}}{q^{m}}}_{m_1 = m_2 = m}
        \underbrace{-\ 2\sum_{m \geq 1}\frac{1}{q^{m+1}}}_{\abs{m_1-m_2}=2} \notag\\
        &\ =\ \frac{2}{q-1}-\frac1q-\frac1{q^2}-\frac{2}{q(q-1)}\ =\ \frac{q-1}{q^2},
    \end{align}
    we get $\Gamma^{\gen}_q = \hat\phi_1(0)\hat\phi_2(0)(q-1)q^{-2}\log^2 (q)/\log^2 (R)+O\paren{\log^{-4}(R)}$. For $p \mid N$ with $p \notin Q_0$ we have $p \gg N^\delta$ and $\abs{G(p)} \leq 2/p$, so $\Gamma^{\gen}_p$ is negligible.

By \eqref{eq:local-cov-values}, $\Gamma_p = 0$ if $p^2 \mid N$, while for $p \Vert N$ only pairs with $m_1,m_2$ both odd contribute, and then $\gamma_{m_1,m_2,N}(p)p^{-(m_1+m_2)/2} = p^{-(m_1+m_2)}$. Replacing $\hat\phi_i$ by $\hat\phi_i(0)$ again costs $O\paren{\log^{-4}(R)}$, and
    \begin{align}
        \sum_{\substack{m_1,m_2 \geq 1\\ m_1,m_2 \text{ odd}}}\frac{1}{q^{m_1+m_2}}\ =\ \left(\sum_{\substack{m \geq 1\\ m \text{ odd}}}\frac{1}{q^{m}}\right)^{\!2}\ =\ \left(\frac{q}{q^2-1}\right)^{\!2},
    \label{eq:cov-value-2}
\end{align}
    so $\Gamma_q = \hat\phi_1(0)\hat\phi_2(0)q^2(q^2-1)^{-2}\log^2 (q)/\log^2 (R)+O\paren{\log^{-4}(R)}$ for $q \in Q$, and $\Gamma_p$ is negligible for $p \Vert N$ with $p \notin Q_0$.

    Combining \eqref{eq:cov-two-blocks} and \eqref{eq:cov-generic-part},
    \begin{align}\label{eq:cov-value-3}
        \Cov(\cF_{N};\phi_1,\phi_2) \ =\ \cE_2(R;\phi_1,\phi_2)-\frac{\hat\phi_1(0)\hat\phi_2(0)}{\log^2 (R)}\left[\sum_{q \in Q_0}\frac{(q-1)\log^2 (q)}{q^2}-\sum_{q \in Q}\frac{q^2\log^2 (q)}{\paren{q^2-1}^2}\right]\notag\\ +O\paren{\log^{-4}(R)},
    \end{align}
    and the bracket is $\kappa_2(\cF)$ by \eqref{eq:kappa2-def}.
\end{proof}

%%%%%%%%%%%%%%%%%%%%%%%%%%%%%%%%%%%%%%%%%%%%%%%%%%%%%%%%%%%%%%%%%%%%%%
\subsection{Evaluating the profiles}\label{subsec:profile-values}

The profiles are now evaluated by partial summation against the Chebyshev function $\theta(t) := \sum_{p \leq t}\log( p)$, with $E(t) := \theta(t)-t$. Define
\begin{align}
    \gamma_{PNT1}\ &:=\ \int_1^\infty \frac{E(t)\paren{1-\log (t)}}{t^2}dt \ \approx\ 2.55511, \label{eq:gPNT1}\\
    \gamma_{PNT2}\ &:=\ \int_1^\infty \frac{E(t)\paren{1-2\log (t)}}{t^3}dt\ \approx\ 0.25840, \label{eq:gPNT2}\\
    \gamma_{PNT3}\ &:=\ 1+\int_1^\infty \frac{E(t)}{t^2}dt\ \approx\ -1.33258, \label{eq:gPNT3}\\
    \gamma_{PNT4}\ &:=\ \int_1^\infty \frac{E(t)\paren{\log^2 (t)-2\log (t)}}{t^2}dt\ \approx\ -10.25383. \label{eq:gPNT4}
\end{align}
All four are Mertens constants for the Chebyshev function. Indeed
\begin{align}
    \gamma_{PNT3} = \lim_{x\to\infty}\paren{\sum_{p \leq x}\log (p)/p-\log (x)} = -\gamma-\sum_p\log (p)/(p(p-1)),
\end{align}
and in the same way $-\gamma_{PNT1}$ and $\gamma_{PNT4}$ are the constants in the asymptotics of $\sum_{p \leq x}\log^2 (p)/p$ and $\sum_{p \leq x}\log^3 (p)/p$, while $1-\gamma_{PNT2} = \sum_p\log^2 (p)/p^2$. Appendix \ref{appendix:partial-summation} gives the derivations and the rapidly convergent series used for the numerical values.

\begin{theorem}\label{thm:profile-values}
    With the notation above,
    \begin{align}
        \cE_1(R;\phi)\ &=\ \frac{\phi(0)}{2}+\frac{2\gamma_{PNT3}\hat\phi(0)}{\log (R)}+\frac{4\gamma_{PNT4}\hat\phi''(0)}{\log^3 (R)}+O\paren{\log^{-5}(R)}, \label{eq:E1-value}\\
        \cE_2(R;\phi_1,\phi_2)\ &=\ \int_0^\infty u\hat\phi_1(u)\hat\phi_2(u)du-\frac{\hat\phi_1(0)\hat\phi_2(0)\paren{\gamma_{PNT1}+1-\gamma_{PNT2}}}{\log^2 (R)}+O\paren{\log^{-4}(R)}. \label{eq:E2-value}
    \end{align}
    Consequently
    \begin{align}\label{eq:D1-profile-value}
        \cD_1(R;\phi)\ =\ \hat\phi(0)+\frac{\phi(0)}{2}+\frac{\paren{\gamma_{k,1}+2\gamma_{PNT3}}\hat\phi(0)}{\log (R)}+\frac{\paren{\tfrac12\gamma_{k,2}+4\gamma_{PNT4}}\hat\phi''(0)}{\log^3 (R)}+O\paren{\log^{-5}(R)},
    \end{align}
    and
    \begin{align}\label{eq:D2-profile-value}
        \cD_2(R;\phi_1,\phi_2)\ =\ &\paren{\hat\phi_1(0)+\frac{\phi_1(0)}{2}}\paren{\hat\phi_2(0)+\frac{\phi_2(0)}{2}} \notag\\ &+2\int_{-\infty}^\infty \abs u\hat\phi_1(u)\hat\phi_2(u)du-2\widehat{\phi_1\phi_2}(0)-\frac{\phi_1(0)\phi_2(0)}{2}+O\paren{\log^{-1}(R)}.
    \end{align}
\end{theorem}

\begin{proof}
    Formula \eqref{eq:E1-value} is Lemma \ref{lem:MSwork} applied with $g(u) = \hat\phi(2u)$, for which $g(0) = \hat\phi(0)$, $g''(0) = 4\hat\phi''(0)$ and $\int_0^\infty g(u)du = \tfrac12\int_0^\infty\hat\phi(v)dv = \tfrac14\int_{-\infty}^\infty\hat\phi = \tfrac14\phi(0)$, the last two equalities using that $\hat\phi$ is even and Fourier inversion.

    For \eqref{eq:E2-value}, isolate the term $(m_1,m_2) = (1,1)$ of \eqref{eq:E2-def}, whose coefficient is $\gamma^{\gen}_{1,1} = 1$:
    \begin{align}
        \sum_p\frac{\log^2 (p)}{p\log^2 (R)}\hat\phi_1\paren{\frac{\log (p)}{\log (R)}}\hat\phi_2\paren{\frac{\log (p)}{\log (R)}}\ =\ \int_0^\infty u\hat\phi_1(u)\hat\phi_2(u)du-\frac{\hat\phi_1(0)\hat\phi_2(0)\gamma_{PNT1}}{\log^2 (R)}+O\paren{\log^{-4}(R)}
    \label{eq:profile-values-1}
\end{align}
    by Lemma \ref{lem:MSwork2} applied to the even function $h = \hat\phi_1\hat\phi_2$. Every remaining pair has $m_1+m_2 \geq 4$, so by Lemma \ref{lem:MSwork3} the remaining terms equal
    \begin{align}
        \frac{\hat\phi_1(0)\hat\phi_2(0)}{\log^2 (R)}\sum_p\log^2 (p)\paren{G(p)-\frac1p}+O\paren{\log^{-4}(R)}\ =\ -\frac{\hat\phi_1(0)\hat\phi_2(0)}{\log^2 (R)}\sum_p\frac{\log^2 (p)}{p^2}+O\paren{\log^{-4}(R)},
    \label{eq:profile-values-2}
\end{align}
    using \eqref{eq:Gq}. As $\sum_p \log^2 (p)/p^2 = 1-\gamma_{PNT2}$, this gives \eqref{eq:E2-value}.

    Now \eqref{eq:D1-profile-value} follows from \eqref{eq:D1-profile-def}, \eqref{eq:U-over-logR} and \eqref{eq:E1-value}. Finally, by \eqref{eq:D2-profile-def} and \eqref{eq:D1-profile-value},
    \begin{align}
        \cD_2(R;\phi_1,\phi_2)\ =\ &\paren{\hat\phi_1(0)+\frac{\phi_1(0)}{2}}\paren{\hat\phi_2(0)+\frac{\phi_2(0)}{2}}-2\paren{\widehat{\phi_1\phi_2}(0)+\frac{\phi_1(0)\phi_2(0)}{2}} \notag\\
        &+4\int_0^\infty u\hat\phi_1(u)\hat\phi_2(u)du+\frac{\phi_1(0)\phi_2(0)}{2}+O\paren{\log^{-1}(R)},
    \label{eq:profile-values-3}
\end{align}
    where we used $(\phi_1\phi_2)(0) = \phi_1(0)\phi_2(0)$. Since $\hat\phi_1\hat\phi_2$ is even, $4\int_0^\infty u\hat\phi_1\hat\phi_2\,du = 2\int_{-\infty}^\infty\abs u\hat\phi_1(u)\hat\phi_2(u)du$, and collecting the terms in $\phi_1(0)\phi_2(0)$ gives \eqref{eq:D2-profile-value}.
\end{proof}

\begin{remark}\label{rmk:katz-sarnak}
    The leading terms of \eqref{eq:D1-profile-value} and \eqref{eq:D2-profile-value} are exactly the one- and two-level densities of the orthogonal symmetry group,
    \begin{align}
        \int_{-\infty}^\infty \phi(x)W_{1,O}(x)dx, \qquad W_{1,O}(x)\ =\ 1+\tfrac12\delta_0(x),
    \label{eq:katz-sarnak-1}
\end{align}
    together with its two-level analogue, for which see \cite{ILS, katz1999random, Mil2}. Our families therefore satisfy the Katz--Sarnak prediction. The profiles $\cD_1$ and $\cD_2$ refine that prediction by recording the full dependence on the conductor while forgetting the factorization entirely, and Theorems \ref{thm:D1-universality-breaks} and \ref{thm:D2-universality-breaks} say that this forgetting is not harmless.
\end{remark}

%%%%%%%%%%%%%%%%%%%%%%%%%%%%%%%%%%%%%%%%%%%%%%%%%%%%%%%%%%%%%%%%%%%%%%
\subsection{Breaking universality} \label{subsec:provingBU}

\FirstUniversalityBreaks*

\begin{proof}
    By \eqref{eq:D1explicitformula}, Lemma \ref{lem:U-universal} and Corollary \ref{cor:S1-master},
    \begin{align}\label{eq:D1-expanded}
        D_1(\cF_{N},\phi)\ =\ \frac{U_k(R;\phi)}{\log (R)}+S_1(\cF_{N},\phi)\ =\ \cD_1(R;\phi)-\frac{2\hat\phi(0)\kappa_1(\cF)}{\log (R)}-\frac{4\hat\phi''(0)\kappa_3(\cF)}{\log^3 (R)}+O\paren{\log^{-5}(R)},
    \end{align}
    which is \eqref{eq:D1-master}. Multiplying by $\log (R)$ and letting $N \to \infty$ gives \eqref{eq:D1-invariant}, whose right side is nonzero as soon as $\hat\phi(0) \neq 0$ and $Q_0(\cF) \neq \emptyset$ by \eqref{eq:kappa1-positive}.

    For \eqref{eq:D1-comparison}, let $R$ and $\tilde R$ be the conductors of the two families. By \eqref{eq:logN-logR} the hypothesis $N \asymp \tilde N$ gives $\log (R)-\log (\tilde R) = O(1)$, hence $\log (R)/\log (\tilde R) \to 1$. By \eqref{eq:D1-profile-value} the profile has the shape $\cD_1(R;\phi) = c_0+c_1\log^{-1}(R)+c_3\log^{-3}(R)+O\paren{\log^{-5}(R)}$ with $c_0,c_1,c_3$ depending only on $k$ and $\phi$, and for any such constant $c$,
    \begin{align} \label{eq:conductor-mismatch}
        \abs{\frac{c}{\log^j (R)}-\frac{c}{\log^j (\tilde R)}}\ =\ \abs c\,\frac{\abs{\log^j (R)-\log^j (\tilde R)}}{\log^j (R)\log^j (\tilde R)}\ \lesssim\ \frac{1}{\log^{j+1} (R)} \qquad (j \geq 1),
    \end{align}
    so $\cD_1(R;\phi)-\cD_1(\tilde R;\phi) = O\paren{\log^{-2}(R)}$. Subtracting the two instances of \eqref{eq:D1-expanded} therefore gives
    \begin{align}
        D_1(\cF_{N},\phi)-D_1(\tilde \cF_{\tilde N},\phi)\ =\ -\frac{2\hat\phi(0)\paren{\kappa_1(\cF)-\kappa_1(\tilde \cF)}}{\log (R)}+O\paren{\log^{-2}(R)},
    \label{eq:subsec-provingbu-1}
\end{align}
    and multiplying by $\log (R)$ yields \eqref{eq:D1-comparison}.
\end{proof}

We turn to the second level. The only ingredient of \eqref{eq:D2-compact} not yet treated is the root-number average $\cW(\cF_N) = \angles{\varepsilon_f}_N$.

\begin{lemma}\label{lem:sign-bias}
    If $N$ is squarefree, then $\cW(\cF_{N}) \ll_{k,\eps} N^{-1/4+\eps}$, and in particular $\cW(\cF_{N})$ is negligible.
\end{lemma}

\begin{proof}
    For squarefree $N$, every $f \in H^*_k(N)$ has $\varepsilon_f = i^k\mu(N)\sqrt{N}\lambda_f(N)$ \cite{ILS}, so
    \begin{align}
        \cW(\cF_{N})\ =\ i^k\mu(N)\sqrt{N}\cdot\frac{\Delta^*_{k,N}(1,N)}{W_R(\cF_{N})}.
    \label{eq:sign-bias-1}
\end{align}
    Since $N$ is squarefree, $N_1 = N$, so \eqref{eq:normalized-general-trace-bound} applies with $D = N$ and $a = b = 1$ and gives $\Delta^*_{k,N}(1,N)/W_R(\cF_{N}) \ll_{k,\eps} N^{-3/4+\eps}$. Multiplying by $\sqrt N$ gives the claim.
\end{proof}

\begin{remark}\label{rmk:nonsquarefree-sign}
    For non-squarefree $N$ the sign $\varepsilon_f$ is determined by Atkin--Lehner eigenvalues at the primes $q$ with $q^2 \mid N$, and it is not a single Hecke eigenvalue. We therefore restrict Theorem \ref{thm:D2-universality-breaks} to squarefree levels. Any bound of the form $\cW(\cF_{N}) = o\paren{\log^{-1}(R)}$ would suffice. If instead $\cW(\cF_{N})$ had a nonzero limit the conclusion would only be strengthened, since $\phi_1(0)\phi_2(0)\cW(\cF_{N})/2$ would then be a further factorization-dependent term of order $1$.
\end{remark}

\SecondUniversalityBreaks*

\begin{proof}
    Since $\hat\phi_1,\hat\phi_2$ are supported in $[-\sigma,\sigma]$ with $\sigma < 2/3$, the function $\phi_1\phi_2$ is even and Schwartz with $\widehat{\phi_1\phi_2} = \hat\phi_1 * \hat\phi_2$ supported in $[-2\sigma,2\sigma]$ and $2\sigma < 4/3$. Hence \eqref{eq:D1-expanded} applies to each of $\phi_1$, $\phi_2$ and $\phi_1\phi_2$. Write
    \begin{align}
        D_1(\cF_{N},\phi)\ =\ \cD_1(R;\phi)+d(\phi), \qquad d(\phi)\ =\ -\frac{2\hat\phi(0)\kappa_1(\cF)}{\log (R)}+O\paren{\log^{-3}(R)}.
    \label{eq:nonsquarefree-sign-1}
\end{align}
    Substituting into \eqref{eq:D2-compact} and using Theorem \ref{thm:cov-value} and Lemma \ref{lem:sign-bias},
    \begin{align}
        D_2(\cF_{N},\phi_1,\phi_2)\ =\ &\cD_1(R;\phi_1)\cD_1(R;\phi_2)+\cD_1(R;\phi_1)d(\phi_2)+\cD_1(R;\phi_2)d(\phi_1)+d(\phi_1)d(\phi_2)\notag\\
        &-2\cD_1(R;\phi_1\phi_2)-2d(\phi_1\phi_2)+4\cE_2(R;\phi_1,\phi_2)+\frac{\phi_1(0)\phi_2(0)}{2}+O\paren{\log^{-2}(R)}\notag\\
        \ =\ &\cD_2(R;\phi_1,\phi_2)+\cD_1(R;\phi_1)d(\phi_2)+\cD_1(R;\phi_2)d(\phi_1)-2d(\phi_1\phi_2)+O\paren{\log^{-2}(R)},
    \label{eq:nonsquarefree-sign-2}
\end{align}
    since $d(\phi_1)d(\phi_2) = O\paren{\log^{-2}(R)}$ and the $\kappa_2$ term of \eqref{eq:cov-value} is $O\paren{\log^{-2}(R)}$. By \eqref{eq:D1-profile-value}, $\cD_1(R;\phi_i) = \hat\phi_i(0)+\tfrac12\phi_i(0)+O\paren{\log^{-1}(R)}$, so the three surviving terms are
    \begin{align}
        -\frac{2\kappa_1(\cF)}{\log (R)}\left[\hat\phi_2(0)\paren{\hat\phi_1(0)+\frac{\phi_1(0)}{2}}+\hat\phi_1(0)\paren{\hat\phi_2(0)+\frac{\phi_2(0)}{2}}-2\widehat{\phi_1\phi_2}(0)\right]+O\paren{\log^{-2}(R)}.
    \label{eq:nonsquarefree-sign-3}
\end{align}
    The bracket equals $2\paren{\hat\phi_1(0)\hat\phi_2(0)-\widehat{\phi_1\phi_2}(0)}+\tfrac12\Upsilon(\phi_1,\phi_2) = 2\Xi(\phi_1,\phi_2)+\tfrac12\Upsilon(\phi_1,\phi_2)$, so the whole expression is $-\kappa_1(\cF)\Xi_2(\phi_1,\phi_2)/\log (R)$ by \eqref{eq:Xi2-def}. This is \eqref{eq:D2-master}, and \eqref{eq:D2-invariant} follows on multiplying by $\log (R)$.

    For \eqref{eq:D2-comparison}, note that by \eqref{eq:D2-profile-def}, \eqref{eq:D1-profile-value} and \eqref{eq:E2-value} the profile has the shape $\cD_2(R;\phi_1,\phi_2) = c_0+c_1\log^{-1}(R)+O\paren{\log^{-2}(R)}$ with $c_0,c_1$ depending only on $k,\phi_1,\phi_2$, so \eqref{eq:conductor-mismatch} gives $\cD_2(R)-\cD_2(\tilde R) = O\paren{\log^{-2}(R)}$ when $N \asymp \tilde N$. Subtracting the two instances of \eqref{eq:D2-master} and multiplying by $\log (R)$ gives \eqref{eq:D2-comparison}.
\end{proof}

%%%%%%%%%%%%%%%%%%%%%%%%%%%%%%%%%%%%%%%%%%%%%%%%%%%%%%%%%%%
\section{Future Work} \label{sec:future-work}

Theorem \ref{thm:D2-universality-breaks} is restricted to squarefree levels. The obstruction is the root number. For squarefree $N$, the sign $\varepsilon_f$ is a Hecke eigenvalue and Lemma \ref{lem:sign-bias} controls its average, whereas for $q^2 \mid N$ the sign is governed by Atkin--Lehner pseudo-eigenvalues, which the trace formula we use does not reach.

Additionally, we require every unbounded prime factor of $N$ to be at least a fixed power of $N$. A prime factor growing more slowly, say of size $\log (N)$, contributes $\log\log (N)/\log (N)$ to the analogue of $\kappa_1(\cF)$, which vanishes in the limit but not fast enough to be absorbed into an error term of size $\log^{-2}(R)$. Such levels should produce a genuinely intermediate behavior, and we do not treat them.

%%%%%%%%%%%%%%%%%%%%%%%%%%%%%%%%%%%%%%%%%%%%%%%%%%%%%%%%%
\appendix

\section{Partial summation estimates}\label{appendix:partial-summation}

Throughout this appendix $\theta(t) = \sum_{p \leq t}\log (p)$ and $E(t) = \theta(t)-t$, and we use only the Prime Number Theorem in the form $E(t) \ll t\exp\paren{-c\sqrt{\log( t)}}$, which guarantees that
\begin{align}\label{eq:E-integrals-converge}
    \int_1^\infty \frac{\abs{E(t)}\log^j (t)}{t^2}dt\ <\ \infty \qquad (0 \leq j \leq 4).
\end{align}
Note also $E(1^-) = -1$. The constants $\gamma_{PNT1},\dots,\gamma_{PNT4}$ are those of \eqref{eq:gPNT1}--\eqref{eq:gPNT4}.

\begin{lemma}\label{lem:MSwork}
    Let $g \colon \bR \to \bR$ be even, four times continuously differentiable, and compactly supported. Then
    \begin{align}\label{eq:MSwork}
        \sum_p \frac{\log (p)}{p}\,g\!\paren{\frac{\log (p)}{\log (R)}}\ =\ \log (R)\int_0^\infty g(u)du+g(0)\gamma_{PNT3}+\frac{g''(0)\gamma_{PNT4}}{2\log^2 (R)}+O\paren{\log^{-4}(R)}.
    \end{align}
    In particular, taking $g(u) = \hat\phi(2u)$ for an even Schwartz $\phi$ with $\hat\phi$ compactly supported,
    \begin{align}\label{eq:E1-lemma}
        \cE_1(R;\phi)\ =\ \frac{\phi(0)}{2}+\frac{2\gamma_{PNT3}\hat\phi(0)}{\log (R)}+\frac{4\gamma_{PNT4}\hat\phi''(0)}{\log^3 (R)}+O\paren{\log^{-5}(R)}.
    \end{align}
\end{lemma}

\begin{proof}
    Write $v = v(t) := \log (t)/\log (R)$ and $h(t) := g(v)/t$, so that the sum is $\int_{1^-}^\infty h(t)\,d\theta(t)$. Since $\theta(t) = t+E(t)$,
    \begin{align}
        \sum_p \frac{\log (p)}{p}g(v(p))\ =\ \int_1^\infty h(t)dt+\int_{1^-}^\infty h(t)\,dE(t).
    \label{eq:mswork-1}
\end{align}
    Substituting $u = v(t)$ in the first integral gives $\log (R)\int_0^\infty g(u)du$. For the second, integrating by parts and using that $h$ has compact support,
    \begin{align}
        \int_{1^-}^\infty h(t)dE(t)\ =\ -h(1)E(1^-)-\int_1^\infty E(t)h'(t)dt\ =\ g(0)-\int_1^\infty E(t)h'(t)dt.
    \label{eq:mswork-2}
\end{align}
    Now $h'(t) = \dfrac{g'(v)}{t^2\log (R)}-\dfrac{g(v)}{t^2}$, so
    \begin{align}
        -\int_1^\infty E(t)h'(t)dt\ =\ \int_1^\infty \frac{E(t)g(v)}{t^2}dt-\frac{1}{\log (R)}\int_1^\infty \frac{E(t)g'(v)}{t^2}dt.
    \label{eq:mswork-3}
\end{align}
    Since $g$ is even, $g'(0) = 0$ and Taylor's theorem gives $g(v) = g(0)+\tfrac12 g''(0)v^2+O(v^4)$ and $g'(v) = g''(0)v+O(v^3)$, with implied constants depending only on $\norm{g^{(4)}}_\infty$. Substituting and using \eqref{eq:E-integrals-converge},
    \begin{align}
        -\int_1^\infty E(t)h'(t)dt\ =\ g(0)\int_1^\infty\frac{E(t)}{t^2}dt+\frac{g''(0)}{2\log^2 (R)}\int_1^\infty\frac{E(t)\paren{\log^2 (t)-2\log (t)}}{t^2}dt+O\paren{\log^{-4}(R)}.
    \label{eq:mswork-4}
\end{align}
    Adding $g(0)$ and recalling \eqref{eq:gPNT3}, \eqref{eq:gPNT4} gives \eqref{eq:MSwork}. Note that only even powers of $1/\log (R)$ occur, because $g$ is even.

    For \eqref{eq:E1-lemma}, take $g(u) = \hat\phi(2u)$, so $g(0) = \hat\phi(0)$, $g''(0) = 4\hat\phi''(0)$, and
    \begin{align}
        \int_0^\infty g(u)du\ =\ \frac12\int_0^\infty \hat\phi(w)dw\ =\ \frac14\int_{-\infty}^\infty \hat\phi(w)dw\ =\ \frac{\phi(0)}{4},
    \label{eq:mswork-5}
\end{align}
    using that $\hat\phi$ is even and Fourier inversion. Multiplying \eqref{eq:MSwork} by $2/\log (R)$ gives \eqref{eq:E1-lemma}.
\end{proof}

\begin{lemma}\label{lem:MSwork2}
    Let $h \colon \bR \to \bC$ be even, twice continuously differentiable with bounded second derivative, and compactly supported. Then
    \begin{align}\label{eq:MSwork2}
        \sum_p \frac{\log^2 (p)}{p}\,h\!\paren{\frac{\log (p)}{\log (R)}}\ =\ \log^2 (R)\int_0^\infty u\,h(u)du-h(0)\gamma_{PNT1}+O\paren{\log^{-2}(R)}.
    \end{align}
\end{lemma}

\begin{proof}
    As above put $v = \log (t)/\log (R)$ and now $h_2(t) := \log( t)\,h(v)/t$, so that the sum is $\int_{1^-}^\infty h_2(t)d\theta(t)$. The main integral is
    \begin{align}
        \int_1^\infty \frac{\log (t)}{t}h(v)dt\ =\ \log^2 (R)\int_0^\infty u\,h(u)du.
    \label{eq:mswork2-1}
\end{align}
    Since $h_2(1) = 0$, integration by parts leaves $-\int_1^\infty E(t)h_2'(t)dt$, and
    \begin{align}
        h_2'(t)\ =\ \frac{h(v)\paren{1-\log (t)}}{t^2}+\frac{\log (t)\,h'(v)}{t^2\log (R)}.
    \label{eq:mswork2-2}
\end{align}
    As $h$ is even, $h(v) = h(0)+O(v^2)$ and $h'(v) = O(v)$, so by \eqref{eq:E-integrals-converge}
    \begin{align}
        -\int_1^\infty E(t)h_2'(t)dt\ =\ -h(0)\int_1^\infty \frac{E(t)\paren{1-\log (t)}}{t^2}dt+O\paren{\log^{-2}(R)},
    \label{eq:mswork2-3}
\end{align}
    which is $-h(0)\gamma_{PNT1}+O\paren{\log^{-2}(R)}$ by \eqref{eq:gPNT1}.
\end{proof}

\begin{lemma}\label{lem:MSwork3}
    Let $\phi_1,\phi_2$ be even Schwartz functions with compactly supported Fourier transforms, and let $(b_{m_1,m_2})_{m_1,m_2 \geq 1}$ be complex numbers vanishing unless $m_1+m_2 \geq 4$ and satisfying
    \begin{align}\label{eq:b-decay}
        \sum_{m_1,m_2 \geq 1}\abs{b_{m_1,m_2}}\paren{1+m_1^2+m_2^2}2^{-(m_1+m_2)/2}\ <\ \infty .
    \end{align}
    Put $B(p) := \sum_{m_1,m_2 \geq 1}b_{m_1,m_2}p^{-(m_1+m_2)/2}$. Then
    \begin{align}
        \sum_p\frac{\log^2 (p)}{\log^2 (R)}\sum_{m_1,m_2 \geq 1}\frac{b_{m_1,m_2}}{p^{(m_1+m_2)/2}}\hat\phi_1\!\paren{\frac{m_1\log (p)}{\log (R)}}\hat\phi_2\!\paren{\frac{m_2\log (p)}{\log (R)}} \notag \\
        =\ \frac{\hat\phi_1(0)\hat\phi_2(0)}{\log^2 (R)}\sum_p \log^2 (p)B(p)+O\paren{\log^{-4}(R)}.
    \label{eq:mswork3-1}
\end{align}
\end{lemma}

\begin{proof}
    Since $\hat\phi_1,\hat\phi_2$ are even, \eqref{eq:phihat-taylor-4} gives
    \begin{align}
        \hat\phi_1\!\paren{\frac{m_1\log (p)}{\log (R)}}\hat\phi_2\!\paren{\frac{m_2\log (p)}{\log( R)}}\ =\ \hat\phi_1(0)\hat\phi_2(0)+O\paren{\paren{m_1^2+m_2^2}\frac{\log^2 (p)}{\log^2 (R)}},
    \label{eq:mswork3-2}
\end{align}
    with an implied constant depending only on $\phi_1,\phi_2$. The resulting error is
    \begin{align}
        \ \ll\ \frac{1}{\log^4 (R)}\sum_p \log^4 (p)\sum_{m_1,m_2 \geq 1}\frac{\abs{b_{m_1,m_2}}\paren{m_1^2+m_2^2}}{p^{(m_1+m_2)/2}}\ \ll\ \frac{1}{\log^4 (R)}\sum_p \frac{\log^4 (p)}{p^{2}},
    \label{eq:mswork3-3}
\end{align}
    where we used $m_1+m_2 \geq 4$ together with $p^{-(m_1+m_2)/2} \leq p^{-2}\,2^{-(m_1+m_2-4)/2}$ and \eqref{eq:b-decay}. The last sum converges, and the main term is as stated. It converges absolutely because $\abs{B(p)} \ll p^{-2}$ by the same estimate.
\end{proof}

\subsection*{Numerical evaluation of the constants}

The integrals \eqref{eq:gPNT1}--\eqref{eq:gPNT4} converge slowly when computed directly, since $E(t)/t$ decays more slowly than any power. They are more efficiently evaluated through the associated Mertens constants, which we record here. Let $\gamma$ be the Euler--Mascheroni constant and $\gamma_1^{\mathrm{St}},\gamma_2^{\mathrm{St}}$ the Stieltjes constants, defined by $\zeta(s) = \tfrac1{s-1}+\gamma-\gamma_1^{\mathrm{St}}(s-1)+\tfrac12\gamma_2^{\mathrm{St}}(s-1)^2-\cdots$. Partial summation as in Lemma \ref{lem:MSwork} gives, for $j = 1,2,3$,
\begin{align}
    \sum_{p \leq x}\frac{\log^{j} (p)}{p}\ =\ \frac{\log^{j} (x)}{j}+c_j+o(1), \qquad
    c_1 = \gamma_{PNT3},\quad c_2 = -\gamma_{PNT1},\quad c_3 = \gamma_{PNT4},
\label{eq:mswork3-4}
\end{align}
and comparing $\sum_{p \leq x}\log^{j}(p)/p$ with $\sum_{n \leq x}\Lambda(n)\log^{j-1}(n)/n$, whose constant terms are read off from the Laurent expansion of $-\zeta'/\zeta$ at $s = 1$, yields the rapidly convergent series
\begin{align}
    \gamma_{PNT3}\ &=\ -\gamma-\sum_p \frac{\log (p)}{p(p-1)}\ =\ -1.3325823\dots,\\
    \gamma_{PNT1}\ &=\ 2\gamma_1^{\mathrm{St}}+\gamma^2+\sum_p \frac{(2p-1)\log^2 (p)}{p(p-1)^2}\ =\ 2.5551076\dots,\\
    \gamma_{PNT4}\ &=\ -2\paren{\tfrac32\gamma_2^{\mathrm{St}}+3\gamma\gamma_1^{\mathrm{St}}+\gamma^3}-\sum_p \log^3 (p)\paren{\frac{p(p+1)}{(p-1)^3}-\frac1p}\ =\ -10.2538271\dots,
\label{eq:mswork3-5}
\end{align}
while partial summation applied to $\sum_p \log^2 (p)/p^2$ gives $\sum_p \log^2 (p)/p^2 = 1-\gamma_{PNT2}$, so that $\gamma_{PNT2} \approx 0.2584021$. We computed the values above from these series using all primes below $5\times10^7$ together with the corresponding integral tail estimates. Each agrees with the direct evaluation of $\sum_{p\le x}\log^j(p)/p-\log^j(x)/j$ at $x = 5\times 10^7$ to within the expected truncation error, and $\gamma_{PNT3}$ agrees with the tabulated value $-1.3325822757$.

%%%%%%%%%%%%%%%%%%%%%%%%%%%%%%%%%%%%%%%%%%%%%%%%%%%%%%%%%%%%%%%%%%%%%%
\section*{Acknowledgements}

This research was conducted at the 2025 and 2026 SMALL REU, with support from the National Science Foundation (DMS-2241623), Columbia University, the Finnerty Fund, Missouri State University, the University of Chicago, the University of Michigan, the University of Rochester, the University of Washington, Williams College, and Yale University.

%%%%%%%%%%%%%%%%%%%%%%%%%%%%%%%%%%%%%%%%%%%%%%%%%%%%%%%%%%%%%%%%%%%%%%
\subsection*{Use of artificial intelligence}

A large language model was used to audit the manuscript for errors, to reorganize and draft the exposition of Sections \ref{chp: S1 compute chpt}--\ref{chp: lower order terms} and Appendix \ref{appendix:partial-summation}, and to carry out the symbolic and numerical verifications reported in 
% Remark \ref{rmk:M30}, 
Theorem \ref{thm:profile-values} and Appendix \ref{appendix:partial-summation}. It produced no mathematical result that the authors did not independently check. All statements, proofs and numerical constants in this paper have been verified by the authors, who take full responsibility for its contents.

%%%%%%%%%%%%%%%%%%%%%%%%%%%%%%%%%%%%%%%%%%%%%%%%%%%%%%%%%

\end{document}